\documentstyle{article}

\newtheorem{thm}{Theorem}[section]
\newtheorem{corol}[thm]{Corollary}
\newtheorem{prop}[thm]{Proposition}

\newtheorem{rem}[thm]{Remark}
\newtheorem{lem}[thm]{Lemma}

\newtheorem{defin}[thm]{Definition}

\newcommand{\beq}{\begin{equation}}
\newcommand{\eeq}{\end{equation}}

\def\NN{\hbox{\sf I\kern-.13em\hbox{N}}}
 \def\HH{\hbox{\sf I\kern-.13em\hbox{H}}}
 \def\DD{\hbox{\sf I\kern-.13em\hbox{D}}}
 \def\RR{\hbox{\sf I\kern-.14em\hbox{R}}}
 \def\CC{\hbox{\sf I\kern-.44em\hbox{C}}}
 \def\ZZ{{\hbox{\sf Z\kern-.43emZ}}}
 \def\QQ{\hbox{\sf C\kern -.48emQ}}
 \def\Cc{\hbox{\sf C\kern -.47em {\raise .48ex \hbox{$\scriptscriptstyle |$}}
   \kern-.5em {\raise .48ex \hbox{$\scriptscriptstyle |$}} }}
 \def\Qq{\hbox{\sf Q\kern -.57em {\raise .48ex \hbox{$\scriptscriptstyle |$}}
   \kern-.55em {\raise .48ex \hbox{$\scriptscriptstyle |$}} }}
\def\AAA{{{\bf A}}}

\title{Seiberg-Witten Floer Homology and Heegaard splittings}
\author{Matilde Marcolli}
\date{}

\begin{document}

\maketitle

{\em Ce fut d'abord une \'etude}\ldots 

{\em Je finis par trouver sacr\'e le d\'esordre de mon esprit.}

(Rimbaud)

\begin{abstract}
\noindent The dimensional reduction of Seiberg-Witten theory defines a
gauge theory of compact connected three-manifolds. Solutions of the
equations modulo gauge symmetries on a three-manifold $Y$ can be
interpreted as the critical points of a functional defined on an
infinite dimensional configuration space of $U(1)$-connections and
spinors. The original Seiberg-Witten equations on the infinite cylinder
$Y\times \RR$ are the downward gradient flow of the functional. Thus,
it is possible to construct an infinite dimensional Morse theory. The
associated Morse homology is the analogue in the context of
Seiberg-Witten theory of Floer's instanton homology constructed using
Yang-Mills gauge theory. The construction and the properties of this
Seiberg-Witten Floer homology are essentially different according to
whether the three-manifold $Y$ is a homology sphere or has non-trivial
rational homology. In this work we construct the Seiberg-Witten Floer
homology for three-manifolds with $b^1(Y)>0$. We define an associated
Casson-like invariant and we prove that it satisfies the expected
intersection formula under a Heegaard splitting of the three-manifold. 
\end{abstract}       

\tableofcontents

\section{Introduction}

The Seiberg-Witten gauge
theory of four-manifold has undergone a rapid and rich development in
the two years following its original formulation by Witten in
1994. Following the  
general strategy illustrated by Atiyah \cite{A}, one is lead to
investigate the relation between gauge theory of three- and of
four-manifolds. 

The Seiberg-Witten gauge theory on four-manifolds has a dimensional
reduction that leads to equations on a three-manifold. As in the case
of Donaldson theory, these equations are the gradient flow of a
functional defined on a Banach manifold. The functional was originally
introduced by Kronheimer and Mrowka in the proof of the Thom
conjecture \cite{KM}. As pointed out by Donaldson \cite{Do2}, it is
possible to define a Seiberg-Witten-Floer homology of a three-manifold by
looking at the critical points and the gradient flow of this
functional. 

For a period of time after the initial introduction of Seiberg-Witten
gauge theory it seemed to be generally accepted that
Seiberg-Witten-Floer homology existed and had the expected properties. 
However, to the best of our knowledge, no detailed
exposition of the construction has yet appeared in the
literature. The purpose of this work was to present in complete detail, the
construction and the properties of the Seiberg-Witten-Floer homology on a 
three-manifold with non-trivial homology 
and of the invariant obtained by computing its Euler
characteristic. A version of this work appeared as a paper
\cite{Ma}. However, the published version is rather sketchy in many
parts, and most of the analytical aspects are not
considered. Moreover, it contains some imprecisions and a mistake that
was pointed out later. The present version corrects the mistake and
fills in the necessary details of the construction.

More recently, a number of results appeared that
enlightened the interesting properties of the Seiberg-Witten Floer
theory and of the associated Casson-like invariant, with important
contributions like \cite{Chen}, \cite{Lim}, \cite{MT},
\cite{MST}, \cite{MOY}, \cite{Wa}, \cite{GWa}.  
Interesting applications have been considered, as in
\cite{Au}, \cite{Au2}, \cite{Mu}. 
Applications to contact structures have been
analysed in \cite{CaWa2}, \cite{KM3}, \cite{LiMa}.
Interesting conjectures arise from the Physics
literature (see e.g. \cite{BL}), where a quantum field theoretic
description of the Floer homology is used. Computations with the three
dimensional equations on a Seifert manifold have been worked out in
\cite{Nico}. 
The study of Seiberg-Witten theory on three-manifolds is a promising
field in rapid expansion and many interesting aspects, like surgery
formulae, are still to be developed.  

\vspace{.3cm}

\noindent 1991 Mathematics Subject Classification: 57R57, 55N35, 58E05.

\subsection*{Acknowledgments} 

\noindent I am grateful to Tom Mrowka for the many
useful conversations, for his very helpful comments and
criticism, and for his patience.  
I am grateful to Bai-Ling Wang for the invaluable and always
stimulating collaboration.  
I thank Rong Wang for having informed me of his interesting work on
the subject, for having pointed out a mistake in a previous version of
this work, and for some very helpful suggestions.
I am grateful to my advisor Mel Rothenberg for having supervised this
work. I thank Stefan Bauer, Rogier Brussee, Ralph Cohen, Rafe Mazzeo,
and Dan Pollack for useful and stimulating conversations.  

\section{Preliminary definitions and notation}

Given a compact connected, oriented 4-manifold $X$,
the Seiberg--Witten equations 
are given in terms of a $U(1)$-connection $\AAA$ and a section $\Psi$ of
the positive spinor bundle $S^+\otimes L$ as

\begin{equation}
D_{\AAA}\Psi=0,
\label{SW1}
\end{equation}
\begin{equation}
F^+_{\AAA}=\Psi\cdot \bar\Psi,
\label{SW2}
\end{equation}
where $D_{\AAA}$ is the Dirac operator twisted with the connection $\AAA$ and
$F_{\AAA}^+$ is the self-dual part of the curvature of $\AAA$.
In the second equation $\Psi\cdot \bar\Psi$ represents the 2-form 
given in local coordinates by

\begin{equation}
(\Psi\cdot\bar\Psi)_{ij}=\frac{1}{4}<e_ie_j\Psi,\Psi>e^i\wedge e^j,
\label{SW2bis}
\end{equation}

\noindent where $<,>$ is the inner product of sections of $S^+\otimes L$. 
The $e_i$
form a local orthonormal basis of sections of $TX$, which acts on $\Psi$ via 
Clifford multiplication. The $e^i$ are the dual basis elements of $T^*X$.

The same equations can be defined on a non-compact 4-manifold $X$ of the form
$X=Y\times \RR$, where 
$Y$ is a compact, connected, oriented 3-manifold without boundary. 

Any oriented three manifold admits a $Spin$-structure, \cite{LaMi}
2.2.3. A choice of the metric determines a natural ``trivial''
$Spin$-structure with spinor bundle $S$. A $Spin_c$-structure is
therefore obtained by twisting $S$ with a line bundle $L$.

Assume a $Spin_c$ structure on $Y$ is given, with spinor bundle 
$S\otimes L$. We can endow the manifold $X$ with the
$Spin_c$ structure defined by
\[ S^\pm\otimes L=\pi_1^* (S\otimes L), \]
where $\pi_1:Y\times \RR\to Y$ is the projection $(y,t)\mapsto y$.

The positive and negative spinor bundles $S^\pm\otimes L$ are isomorphic
via Clifford multiplication by $dt$. 

\begin{defin}
A pair $(\AAA,\Psi)$ is {\em in a temporal gauge} if the $dt$ component of
the connection $\AAA$ is identically zero.
\label{tempgauge}
\end{defin}

\begin{rem}
Any pair $(\AAA,\Psi)$ is gauge equivalent (on $X$) to a pair in a
temporal gauge. 
\label{r1}
\end{rem}

An element $(\AAA,\Psi)$ in a
temporal gauge on $Y\times \RR$ can therefore be written as a path
$(A(t),\psi(t))$.
We obtain the dimensional reduction of the gauge theory as follows.

\begin{lem}
For a pair $(\AAA,\Psi)$ in a temporal gauge,
the Seiberg--Witten equations (\ref{SW1}) and (\ref{SW2}) induce the following
equations on $Y$:
\begin{equation}
\frac{d}{dt}\psi(t)=-\partial_{A(t)}\psi(t),
\label{3SW1}
\end{equation}
and
\begin{equation}
\frac{d}{dt}A(t)=-*F_{A(t)}+\sigma(\psi(t),\psi(t)),
\label{3SW2}
\end{equation}
where the imaginary 1-form $\sigma(\psi,\psi)$ is given in local 
coordinates by $\frac{1}{2}<e_i\psi,\psi>e^i$.
\label{dimred}
\end{lem}

\noindent\underline{Proof:} The first equation (\ref{3SW1}) is obtained
by writing the Dirac operator on $X$ as $D_A=\partial_t
+\partial_{A(t)}$, where $\partial_{A(t)}$ is the self-adjoint Dirac
operator on $Y$ twisted with a time dependent connection $A(t)$.

To obtain the equation (\ref{3SW2}), write the equation
(\ref{SW2}) in local coordinates and consider separately the 
basis elements that are in $\Lambda^2 T^*Y$ and those of the form
$e^i\wedge dt$, with $e^i$ the local basis of $T^*Y$.

This gives equations of the form
\[ \frac{1}{2}(F_{it}+\epsilon^{itjk} F_{jk}) e^i\wedge dt
=\frac{1}{4}<e_ie_t\psi,\psi> e^i\wedge dt, \] 
where we used that $F^+=\frac{1}{2}(F+*F)$ on the 4-manifold $X$.
The symbol $\epsilon$ is the sign of the permutation $\{ itjk \}$ in
$\Sigma_4$. 
Up to composing with the $*$-isomorphism on $Y$ and identifying the positive
and negative spinors via Clifford multiplication by $dt$, these are
the equations (\ref{3SW2}) of the reduced gauge theory. 

\noindent QED 

The gauge group ${\cal G}$ acting on the space of solutions 
on $Y$ will be the subgroup of ${\cal G}(X)$, the group 
of gauge transformations on $X$, which consists of translation invariant
gauge transformations, $\frac{d\lambda}{dt}=0$.

The action of the gauge group is given by
\[ (A,\psi)\mapsto (A-i\lambda^{-1}d\lambda, \lambda\psi). \]

The equations (\ref{3SW1}), (\ref{3SW2}) represent the gradient flow of a 
functional, defined on the space of connections on $Y$ and sections of
$S\otimes L$. This has been introduced in \cite{LM}, more
details have been worked out for instance in \cite{Au}. The Floer
theory of this functional and its properties have been considered
independently by various authors \cite{CaWa}, \cite{Ma}, \cite{Wa},
\cite{GWa}. 

\begin{defin}
For a fixed connection $A_0$, we define a functional
\begin{equation}
{\cal C}(A,\psi)=\frac{-1}{2}\int_Y (A-A_0)\wedge (F_A+F_{A_0})
+\frac{1}{2}\int_Y <\psi,\partial_A\psi>dv.
\label{funct}
\end{equation}
\end{defin}

Upon a gauge transformation the functional (\ref{funct}) changes according
to
\[ {\cal C}(A-i\lambda^{-1}d\lambda, i\lambda\psi)={\cal C}(A,\psi)-2i\pi
\int_Y c_1(L)\wedge \lambda^{-1}d\lambda, \]
where $c_1(L)$ is the representative of the first Chern class given by
the curvature 2-form.

Let $h(\lambda)=[\frac{i}{2\pi}\lambda^{-1}d\lambda]$ be the class in
$H^1(Y,\ZZ)$ representing the the Cartan-Maurer form of
$\lambda$. Notice that $\pi_0({\cal G})\cong H^1(Y,\ZZ)$.
Since the cohomology class $h(\lambda)$
represents the connected component of the gauge transformation $\lambda$
in the gauge group, the functional (\ref{funct}) is well defined on
the space ${\cal B}^0={\cal A}/{\cal G}^0$ of connections and
sections, modulo the action of the trivial connected
component ${\cal G}^0$ of the gauge group. In some cases, however,
this configuration space turns out to be ``too large'', in the sense
that the corresponding moduli space can loose the compactness. This
happens for instance in the case when $c_1(L)=0$, as will be discussed
later. It is therefore convenient to introduce another configuration
space which is larger than ${\cal B}$ but smaller than ${\cal
  B}^0$, on which (\ref{funct}) is well defined as an
$\RR$-valued functional. This was first introduced by R.G. Wang \cite{GWa}.

Consider the subgroup of ${\cal G}$ 
\beq
\tilde{\cal G}=\{ \lambda\in {\cal G} | ~\frac{i}{2\pi}\int_Y
c_1(L)\wedge \lambda^{-1}d\lambda=0. \}. 
\label{int=0}
\eeq
The subgroup $\tilde{\cal G}$ is the kernel of the homomorphism
\[ \xi: {\cal G}\to \ZZ \]
\[ \xi(\lambda)=\frac{i}{2\pi}\int_Y
c_1(L)\wedge \lambda^{-1}d\lambda. \]
The induced homomophism $\xi :H^1(Y,\ZZ)\to\ZZ$ determines the
subgroup
\[ \tilde H^1(Y,\ZZ)\equiv \{ h\in H^1(Y,\ZZ)\mid 
< c_1(L)\cup h,[Y]>=0\}. \]
Let $H=H^1(Y,\ZZ)/ \tilde H^1(Y,\ZZ)$. We consider the configuration
space ${\cal B}_H={\cal A}/{\tilde{\cal G}}$. This is a covering of
${\cal B}$ with covering group $H$.

In particular let $\hat{\cal A}$ be the space of irreducible
pairs $(A,\psi)$, where $\psi$ is not identically zero. We consider
this space modulo the
action of  $\tilde{\cal G}$, 
\[ \hat{\cal B}_H=\hat{\cal A}/{\tilde{\cal G}}, \]
and restrict the functional ${\cal C}$ of (\ref{funct}) over $\hat
{\cal B}_H$. 

We choose to topologize the space of connections and sections
$\hat{\cal A}$ with a fixed $L^2_k$ Sobolev norm, with $k>2$, and 
the space of gauge transformations ${\cal G}$ with the $L^2_{k+1}$
norm. This makes $\hat{\cal B}_H$ into a Banach manifold.

The functional (\ref{funct}) is the analogue of the
Chern--Simons functional in Donaldson theory.

\section{Extremals and Morse Index}

The critical points of the functional ${\cal C}$ of (\ref{funct}) are
the pairs $(A, \psi)$ up to gauge transformation that satisfy
\begin{equation}
\begin{array}{c}
         \partial_A\psi=0, \\
         {}*{}F_A=\sigma(\psi,\psi).
\end{array}
\label{extrem}
\end{equation}

In this section we prove that under a suitable perturbation the space
of critical points of the functional ${\cal C}$ modulo  
gauge transformations is an oriented compact zero dimensional manifold.
We compute the Hessian of the functional and
define the relative Morse index of critical points.

\subsection{Deformation Complex}

Let ${\cal M}_C$ denote the set of critical points of ${\cal C}$ in
the Banach manifold $\hat{\cal B}$, i.e. the set
of solutions of (\ref{extrem}) modulo gauge. In the following we shall
always write $\Lambda^*$ for the complex of imaginary valued forms.
The virtual tangent space of ${\cal M}_C$ is given by 
\[ Ker(T)/Im(G), \]
where the operator 
\[ G: \Lambda^0_{L^2_{k+1}}(Y)\to
\Lambda^1_{L^2_k}(Y)\oplus\Gamma_{L^2_k} (S\otimes L), \]
\begin{equation}
G\mid_{(A,\psi)}(f)=(-df,f\psi)
\label{G}
\end{equation}
is the infinitesimal action of
the gauge group and the map $T$ is the linearization of the 
equation (\ref{extrem}) at a pair $(A_0,\psi_0)$,
\[ T: \Lambda^1_{L^2_k}(Y)\oplus\Gamma_{L^2_k} (S
\otimes L)\to \Lambda^1_{L^2_{k-1}}(Y)\oplus\Gamma_{L^2_{k-1}} (S
\otimes L), \]
\begin{equation}
T\mid_{(A,\psi)}(\alpha,\phi)=\left\{
\begin{array}{c}
      *d\alpha -\sigma(\psi,\phi)-\sigma(\phi,\psi) \\
      \partial_{A}\phi +\alpha \psi.
\end{array}\right.
\label{linextr}
\end{equation}

These operators fit into an elliptic complex, the {\em deformation
  complex} $C^*$  of equation (\ref{extrem})
\[ 0\to\Lambda^0_{L^2_k}(Y) \oplus \Lambda^1_{L^2_k}(Y) \oplus
\Gamma_{L^2_k} (S\otimes L)\stackrel{L}{\to}
\Lambda^0_{L^2_{k-1}}(Y) \oplus \Lambda^1_{L^2_{k-1}}(Y)\oplus
\Gamma_{L^2_{k-1}} (S\otimes L)\to 0, \]
where the operator $L$ is
\[ L\mid_{(A,\psi)}(f,\alpha,\phi)=\left\{\begin{array}{l}
            T\mid_{(A,\psi)}(\alpha,\phi)+G\mid_{(A,\psi)}(f)\\
            G^*\mid_{(A,\psi)}(\alpha,\phi)
          \end{array}\right. \]

Here the operator $G^*_{A,\psi}$ is the formal adjoint of the
linearization of the group action $G_{A,\psi}$.  We have
\[ G^*\mid_{(A,\psi)}(\alpha,\phi)=-d^*\alpha +i Im <\psi,\phi >. \]

\begin{lem}
The operator $L$ has index zero.
\end{lem}

\noindent\underline{Proof:} Up to zero order operators, $L$ reduces to the
Dirac operator (which has index zero on a three manifold), and the 
elliptic complex
\[ 0\to\Lambda^0(Y)\oplus\Lambda^1(Y)\stackrel{D}{\to}\Lambda^0(Y)
\oplus\Lambda^1(Y)\to 0,\]
where $D$ is given by 
\begin{equation}
D=\left(\begin{array}{cc}
           0&d^*\\
           d&*d
         \end{array}\right). 
\label{D}
\end{equation}
The index of $D$ is the Euler characteristic that is trivial on a
closed three manifold $Y$.

\noindent QED

The virtual dimension of the space of critical points is given by the
first Betti number $h^1(C^*)$ of the short complex. 
Unlike the four dimensional problem, since
\[ 0=Ind(L)=-\chi(C^*)=h^1(C^*)-h^2(C^*), \]
the index computation does not give enough information on the dimension
$h^1(C^*)$. We shall introduce a suitable perturbation so that we can
always restrict to the case with $h^1(C^*)=0$, that is to
a zero-dimensional moduli space.

\subsection{Compactness}

In order to prove that the set of critical points of $\tilde {\cal C}$ modulo
gauge transformations is compact  we use 
Sobolev techniques.

\begin{thm}
Any sequence $\{ (A_j,\psi_j) \}$ of solutions
of (\ref{extrem}) has a subsequence that converges with all
derivatives (up to gauge transformations) to another solution.
\label{compact}
\end{thm}

It is enough to show that there exists a sequence of gauge
transformations $\{ \lambda_j \}$ such that all $L^2_k$ norms of
$(A_j-i\lambda_j^{-1}d\lambda_j, \lambda_j\psi_j)$ are bounded. The
result then follows by the Sobolev embedding theorem.

\begin{lem}
A section $\psi$ of the spinor bundle $S\otimes L$ that is a
solution of (\ref{extrem}) has bounded $L^2$ norm on $Y$.
\label{boundsect}
\end{lem}    

The argument follows the line of the analogous 
result proved by Kronheimer and Mrowka \cite{KM}
for the Seiberg-Witten equations on a compact 4 manifold. 
The inequality
\[ \mid \psi\mid \leq \max_Y (0,-\kappa), \]
holds, with $\kappa$ the scalar curvature.  

\begin{lem}
Suppose given a sequence $\{ (A_j,\psi_j) \}$ of solutions of (\ref{extrem})
on $Y$, and a fixed connection $A_0$. Then there exists a sequence of gauge
transformations $\lambda_j$, 
such that the $A_j-A_0-i\lambda_j^{-1}d\lambda_j$ are
co-closed 1-forms.
\label{1coclosed}
\end{lem}

\noindent\underline{Proof:} Consider the tangent space to the
set of solutions of (\ref{extrem}). At the connection $A_0$ the
subspace spanned by the infinitesimal action of the gauge group is the
image of $d$. 
Thus we can find a sequence of gauge transformations that makes
the 1-forms $A_j-A_0$ orthogonal to the image of $d$, i.e. in
the kernel of $d^*$.

\noindent QED

We also need the following gauge fixing condition (this is the
analogue of the four dimensional gauge fixing condition, \cite{Mo}
Lemma 5.3.1.

\begin{lem}
The gauge transformations $\lambda_i$ of Lemma \ref{1coclosed} can be
chosen so that the sequence $\lambda_i (A_i)-A_0$ satisfies
\[ \| \lambda_iA_i-A_0\|^2_{L^2_k}\leq C\| F_{\lambda_iA_i}
\|_{L^2_{k-1}}+K. \] 
\label{gaugefix}
\end{lem}

\noindent\underline{Proof:} 
We follow the analogous argument given \cite{Mo} in the
four-dimensional case. Consider the operator
\[ (d^*,d): \Lambda^1_{L^2_k}  \to
\Lambda^2_{L^2_{k-1}}\oplus\Lambda^0_{L^2_{k-1}}, \]
and a decomposition of the 1-form $\lambda A -A_0$ in a harmonic
component and a component that is orthogonal to the harmonic forms,
\[ \lambda A -A_0=h+\beta. \]
Then $\beta$ satisfies
\[ \| \beta \|^2_{L^2_k}\leq C\| (d^*(\beta),d(\beta))
\|^2_{L^2_{k-1}}. \]
Since we can assume that $d^*(\beta)=0$ by Lemma \ref{1coclosed}, and
$d(\lambda A -A_0)=F_A-F_{A_0}$, we have
\[ \| \beta \|^2_{L^2_k}\leq C\| F_A \|^2_{L^2_{k-1}}+ C_1. \]
The harmonic component $h$ may not be bounded, but we can rewrite
$h=h_1+h_2$, where $h_2$ is a harmonic one-form in the lattice $\Sigma$ of
imaginary valued harmonic 1-forms with periods in $2\pi i \ZZ$. The
form $h_1$ is in the quotient $H^1(Y,i\RR)/\Sigma$, which is a compact
torus. Thus, the norm of $h_1$ is bounded, $\| h_1 \|\leq C_2$, and
the term $h_2$ can be eliminated with a gauge transformation. In fact
$h_2$ can be written as $h_2=h(\lambda_2)$ for some $\lambda_2: Y\to
U(1)$, and $\lambda_2(\lambda A- A_0)$ satisfies the required
estimate, with constant $K= CC_2+C_1$.

\noindent QED 

Notice that Lemma \ref{gaugefix} requires the use of the full gauge
group ${\cal G}$, and it no longer holds if we restrict to the
identity component ${\cal G}^0$.

\noindent\underline{Proof of theorem \ref{compact}:}
For simplicity of notation we can assume that the forms $A_j-A_0$ are
coclosed and satisfy the estimate of Lemma \ref{gaugefix}. We can
write $(d+d^*)A_j$ instead of $dA_j$ and  
$(dd^*+d^*d)A_j$ instead of $d^*d A_j$. Since $d+d^*$ and $dd^*+d^*d$
are elliptic operators, we can use elliptic estimates to bound Sobolev
norms.

In fact by lemma \ref{boundsect} and the second equation of
(\ref{extrem}) we have a bound on the $L^2$
norm of $dA_j=F_{A_j}$. 
This immediately gives a bound on the $L^2_1$
norm of $A_j$, because of Lemma \ref{gaugefix}. 
The first equation in (\ref{extrem}) and lemma \ref{boundsect}
give a bound on the $L^2_1$ norm of the sections $\psi_j$,
since the Dirac operator is elliptic. In fact we have an
estimate
\[ \| \psi_j \|_{L^2_1}\leq c(\| \partial_A\psi_j \|_{L^2}+\| \psi_j
\|_{L^2}). \]   

Moreover, we can obtain a bound on the $L^2$ norm of $d^*F_{A_j}$. This
can be computed in local coordinates from the second equation of
(\ref{extrem}). 

We have
\[ \| d^*F \|_{L^2}\leq\sum_{i\neq j} \mid 2<\psi,e_i\nabla_j\psi >\mid^2 \]
and 
\[ \mid <\psi,e_i\nabla_j\psi> \mid\leq c \|\psi \|_{L^2_1}
\| \psi \|_{L^2}. \]
Therefore we obtain a bound on $\| d^*F_{A_j} \|_{L^2}$ in terms of
the bounds on $\| \psi \|_{L^2}$ and $\| \psi \|_{L^2_1}$.

Now we can use a bootstrapping argument to bound higher Sobolev norms;
for instance we get
\[  \| A_j \|_{L^2_2}\leq c(\| d^*d A_j \|_{L^2}+ \| A_j \|_{L^2}). \]

We can bound the higher Sobolev norms of $\psi_j$ by the elliptic
estimates applied to the Dirac operator (see the analogous argument
in the four-dimensional case \cite{Mo}). 

\noindent QED 

\begin{corol}
A consequence of theorem \ref{compact} is that we can improve the
degree of regularity of the elements in ${\cal M}_C$. Namely, if we
consider the moduli space ${\cal M}_C^{k'}$ inside $\hat{\cal
  B}_{L^2_{k'}}$ with $k'>k$, then the natural map ${\cal M}_C^{k'}\to
{\cal M}_C^{k}$ is a diffeomorphism.
\label{morgan}
\end{corol} 

The proof of corollary \ref{morgan} follows as in the four-dimensional
case \cite{Mo}.

Theorem \ref{compact} holds when we consider critical points in the
space $\hat{\cal B}$, i.e. when we identify points under the action of
$H^1(Y,\ZZ)$ on $\hat{\cal B}^0$. For manifolds $Y$ with $b^1(Y)>0$,
the moduli space ${\cal M}_C^0$ in the configuration space ${\cal
  B}^0$ is non-compact. In fact, it is an $H^1(Y,\ZZ)$ covering of
${\cal M}_C$. The moduli space ${\cal M}_C^H$ inside ${\cal B}_H$ is
compact iff the quotient group $H$ is finite. In fact, ${\cal M}_C^H$
is an $H$ covering of ${\cal M}_C$. Notice that the quotient $H$ is
either trivial or a copy of $\ZZ$, hence the moduli space ${\cal
  M}_C^H$ is either ${\cal M}_C$ or a $\ZZ$-covering of it, depending
on the Chern class $c_1(L)$.

\subsection{Perturbation}

The condition $Ker(L)=0$ in the deformation complex
ensures surjectivity. Therefore, when it is satisfied at all solutions, 
by the implicit function theorem on Banach manifolds we have that ${\cal M}_C$
is a smoothly embedded submanifold of ${\cal B}$. 
However, in general we also have to consider the presence of 
reducible solutions, i.e. solutions with $\psi\equiv 0$ and a non-trivial
stabiliser of the action of the gauge group (the group of constant gauge
transformations). The corresponding moduli space fails to be a smooth
manifold at such points.
Thus, in order to avoid reducible solutions, we consider a perturbed
functional, so that the linearization at a solution of the perturbed equations
will have $L$ surjective.
There is a natural choice of the perturbation, which is simply the three
dimensional reduction of the perturbed Seiberg--Witten equations in
four dimensions \cite{Au}, \cite{KM}, \cite{W2}.
In four dimensions these are
\[ D_A\psi=0, \]
\[ F^+_A=\psi\cdot\bar\psi +\delta, \]
where $\delta$ is a self-dual imaginary 2-form.

The corresponding dimensional reduction gives
\begin{equation}
\frac{d}{dt}\psi=-\partial_A\psi 
\label{3SW1P}
\end{equation}
and
\begin{equation}
\frac{d}{dt}A=\sigma(\psi,\psi)-*F_A+2\rho, 
\label{3SW2P}
\end{equation}
where $\rho$ is an imaginary 1-form on $Y$ such that
$\delta=dt\wedge\rho +*(dt\wedge \rho)$, as explained in \cite{Au}.

These equations can be thought of as the downward gradient flow of a functional
\begin{equation}
\tilde {\cal C}(A,\psi)={\cal C}(A,\psi)-2\int_Y (A-A_0)\wedge *\rho.
\label{pertfunct}
\end{equation}

The functional (\ref{pertfunct}) has the same behaviour as (\ref{funct}) 
under gauge transformations, since the form $\rho$ is divergence-free.
In fact, it is easy to see that the equations
\[ \partial_A \psi =0, \]
\[ *F_A=\sigma(\psi,\psi)+2i\rho \]
and the expression
\[ d^*\sigma(\psi,\psi)= \frac{1}{2} (<\partial_A\psi,\psi>-
<\psi,\partial_A\psi>)  \]
imply that the perturbation $\rho$ is co-closed.

The critical points of the perturbed functional (\ref{pertfunct})
satisfy the equations
\beq
\begin{array}{c}
\partial_A\psi=0 \\[2mm]
{}*F_A-\sigma(\psi,\psi)-2\rho=0.
\end{array}
\label{extremP}
\eeq

In order to show that this choice of perturbation has the
expected property we need the following simple computation.

\begin{lem}
The following identity holds: 
\[ \frac{1}{2}\int_Y <\alpha \psi,\phi> dv =\int_Y \alpha \wedge *\sigma
(\psi,\phi)=-<\alpha,\sigma(\psi,\phi)>,\]
where $\psi$ and $\phi$ are sections of $S\otimes L$ and
$\alpha$ is a 1-form.
\label{ide}
\end{lem}

\noindent\underline{Proof:} In local coordinates 
\[ *\sigma(\psi,\phi)=\frac{1}{2}<e_i\psi,\phi>e^j\wedge e^k, \]
\[ \alpha\wedge *\sigma(\psi,\phi)=<\alpha_ie_i\psi,\phi>dv. \]

\noindent QED

We use it to prove the following lemma.

\begin{lem}
Assume that the Chern class $c_1(L)$ is non-trivial.
For a generic (Baire second category) set of co-closed perturbations $\rho$
that satisfy the condition $*\rho\neq i\pi c_1(L)$, no 
reducible solution arises among the critical points
of the functional (\ref{pertfunct}).  Moreover, the operator $L$ of the 
deformed complex $\tilde C^*$ obtained by linearization at a solution of
the perturbed equations is surjective.
Hence the dimension of the critical set modulo the action of
the gauge group is zero, since $h^1(\tilde C^*)=0$.
\label{dim1}
\end{lem}

\noindent\underline{Proof:} 
Consider the operator
\[ \tilde L\mid_{(A,\psi,\rho)}(\alpha,\phi,\eta)=-2\eta + L(\alpha,\phi). \]
We first prove that this operator is surjective. We know that $L$ is a 
Fredholm operator, hence $\tilde L$ has a closed range. Therefore it
is sufficient to prove that $\tilde L$ has dense range.

Let $(\beta,\xi,g)$ be an element
that is $L^2$-orthogonal to the range of $\tilde L$. Then
$(\beta,\xi,g)$ is in the kernel of the adjoint, hence by elliptic
regularity we can consider the $L^2$ pairing of $L^2_k$ and $L^2_{-k}$,
\[ \langle \beta, -*d\alpha -df +\sigma(\psi,\phi) +\sigma(\phi,\psi)-
2\eta \rangle + 
\langle \xi, \partial_A \phi+\alpha\psi+f\psi \rangle+ \]
\[ +\langle g,G^*(\alpha,\phi)\rangle=0. \]

By varying $\eta$ we force $\beta\equiv 0$. The vanishing of 
\[ \langle \xi, \partial_A\phi+\alpha\psi+f\psi \rangle+ 
\langle g,G^*(\alpha,\phi)\rangle \]
gives an equation $\Delta g +1/2 g |\psi|^2=0$ which implies $g\equiv
0$ by the maximum principle. Then by varying $\phi$ and $\alpha$ we
get $\partial_A\xi=0$ and $\sigma(\xi,\psi)+\sigma(\psi,\xi)=0$. The
latter is satisfied if $\xi$ is an imaginary multiple of $\psi$,
$\xi=i\lambda\psi$, where neither of the two vanishes. Both $\xi$ and
$\psi$ are in the kernel of $\partial_A$, thus if either of them
vanishes on an open set it has to vanish identically (and we know that
$\psi$ is not identically zero). If we have $\xi=i\lambda\psi$, we
obtain that $\xi$ is identically zero as a consequence of the
vanishing of the inner product $\langle \xi, f\psi\rangle$ for
arbitrary smooth compactly supported functions $f$.

Consider the set ${\cal W}$ that fibres over the co-closed 1-forms 
assigning to 
$\rho$ the set of solutions of the corresponding perturbed equations modulo 
gauge. Since $\tilde L$ is surjective, zero is a regular value and therefore
${\cal W}$ is a smooth manifold.
Now by standard Fredholm techniques we see that the projection
map
\[ {\cal W}\to Ker(d^*)\subset\Lambda^1(Y), \]
\[ (A,\psi,\rho)\mapsto \rho, \]
linearizes to a Fredholm projection
\[ Ker(\tilde L)\to \Lambda^1(Y). \]

This implies that the Morse Sard lemma applies, and the set of regular values 
of the projection map is a Baire second category set. At a regular
value $\rho$ the set of solutions of (\ref{extrem}) modulo gauge is
therefore a smooth submanifold of ${\cal W}$ which is cut out
transversally. The transversality property ensures that the operator
$L$ computed at these solutions is surjective.
 
The condition $*\rho\neq \pi c_1(L)$, which guarantees that no reducible 
solution arises, comes from the analogous result in the
four dimensional theory and still gives an open dense condition.

\noindent QED 

\begin{corol}
The moduli space ${\cal M}_C$ of critical points of $\tilde {\cal C}$
modulo gauge transformations is a discrete set of points. The moduli
spaces ${\cal M}_C^0$ and ${\cal M}_C^H$ are given by
copies of ${\cal M}_C$ for each point of $H^1(Y,\ZZ)$ and $H$
respectively. 
\label{discrete}
\end{corol}

Suppose that $c_1(L)$ vanishes but we still have $b^1(Y)>0$. In this
case  reducible solutions can be
perturbed away adding a harmonic 1-form to the curvature
equation. In the case of a homology sphere reducible solutions can not
be perturbed away. In fact, for a perturbation $\rho=*d\theta$ we
would have the reducible solution $(\theta, 0)$. The case of manifolds
with $b^1(Y)=0$ has been analysed in \cite{Wa}. It is remarkable that
for rational homology spheres the reducible solution leads to a subtle
metric dependence problem which is not present in the case with $b^1(Y)>0$.

\subsection{Orientation}

An orientation of the moduli space of critical points of $\tilde {\cal C}$
can be constructed following the procedure illustrated in \cite{D}.
A trivialization of the
determinant line bundle $Det(L)\mid_{(A,\psi)}$ determines an orientation
of the tangent bundle of ${\cal M}_C$ if we show that the action of the
gauge group is orientation preserving. 
We can orient $Det(\partial_A)$ with the orientation induced on the
bundle of spinors by the 
complex structure of the bundle $S^\pm$ on $X$.
A choice of an orientation for
the vector space $H^1(Y,\RR)$ will orient $Det(D)$, where
$D$ is the operator defined in (\ref{D}).  This gives an orientation
of the determinant line bundle $Det(L)\mid_{(A,\psi)}$. 

In order to have an induced orientation on the moduli space
it is sufficient to show that the action of the gauge group preserves
the orientation.

\begin{lem}
The action of the gauge group on the space of solutions of
(\ref{extrem}) preserves the orientation of the tangent space.
\label{or}
\end{lem}

\noindent\underline{Proof:} The action of ${\cal G}$, $\lambda
(A,\psi)= (A-i\lambda^{-1}d\lambda, \lambda\psi)$, induces an
isomorphism 
\[ \lambda^*: Ker(L\mid_{(A,\psi)})\to
Ker(L\mid_{\lambda(A,\psi)}), \]
\[ (f,\alpha,\phi)\mapsto (f,\alpha,\lambda\phi). \]
This map is complex linear on the tangent space of $S$, 
hence the orientation is well defined on ${\cal M}_C$.

\noindent QED    

\subsection{The Hessian}

The critical points of the functional ${\cal C}$ (or of the perturbed one
$\tilde {\cal C}$), which means the
translation invariant solutions of (\ref{3SW1}) and (\ref{3SW2}) on
$Y\times \RR$, 
play the same role as the flat connections in Donaldson theory. Flat
connections are in
fact the critical points of the Chern--Simons functional.

Some technical difficulties arise in defining the Morse index of a
critical point.  
This is defined as the number of negative eigenvalues
of the Hessian. In order to compute the Hessian we can use the
following results.

\begin{thm}
Consider a parametrised family 
\[ (A_s,\psi_s)=(A,\psi)+s(\alpha,\phi) \] 
of connections and sections. The linear part of the increment, that
is the coefficient of $s$ in
\[ \tilde {\cal C}(A_s,\psi_s)-\tilde {\cal C}(A,\psi), \]
is a 1-form on the infinite dimensional space of connections and
spinors which is identified via the metric with
the gradient of $\tilde {\cal C}$. The coefficient of $s^2/2$ is 
a 2-form which induces an operator $T$ on the tangent space of ${\cal B}$ at
$(A,\psi)$. This operator can be identified via the metric with the
linearization (\ref{linextr}) of the equations of the gauge theory. At
a critical point this will be the Hessian. 
\label{hessian}
\end{thm}

\noindent\underline{Proof:} The explicit form of the increment
$\tilde {\cal C}(A_s,\psi_s)-\tilde {\cal C}(A,\psi)$ is 
\[ \frac{s}{2}(-\int_Y \alpha\wedge F_A + (A-A_0)\wedge d\alpha + \]
\[ \int_Y (<\psi,\alpha\psi> + <\phi, \partial_A\psi> +
<\psi,\partial_A\phi>) dv \]
\[ +4\int_Y \alpha\wedge *\rho ) + \]
\[ \frac{s^2}{2} ( -\int_Y \alpha\wedge d\alpha + \int_Y
(<\phi,\alpha\psi>+ <\psi,\alpha\phi>)dv \]
\[ + \int_Y <\phi,\partial_A\phi>dv), \]
plus higher order terms.

Using lemma \ref{ide} we can write the first order increment as
\begin{equation}
\label{1forma}
\begin{array}{c}
{\cal F}\mid_{(A,\psi,\rho)}(\alpha,\phi)=-\int_Y i\alpha\wedge 
(F_A -2*\rho -*\sigma(\psi,\psi)) + \\[2mm]
+\frac{1}{2}\int_Y <\phi,\partial_A\psi>+<\partial_A\psi,\phi>.  
\end{array}
\end{equation}
This is the inner product of the tangent vector $(\alpha,\phi)$ with
the gradient flow  (\ref{3SW1P}), (\ref{3SW2P}).

The Hessian is a quadratic form in the increment $(\alpha,\phi)$,
which is a vector in the $L^2_k$-tangent space. By
means of lemma \ref{ide} the quadratic term can be rewritten as
\[ \frac{s^2}{2}( -\int_Y \alpha\wedge d\alpha + 2Re\int_Y <\phi,\alpha\psi> +
\int_Y <\phi,\partial_A\phi>) = \]
\[ \frac{s^2}{2}( -\int_Y \alpha\wedge d\alpha +\int_Y\alpha\wedge
(*\sigma(\phi,\psi)+*\sigma(\psi,\phi)  + Re \int_Y
<\phi,\alpha\psi+\partial_A\phi> ). \] 
Thus we have the Hessian of the form
\begin{equation}
\label{2forma}
\nabla{\cal F}\mid_{(A,\psi,\rho)}(\alpha,\phi)=<\alpha, *d\alpha 
-\sigma(\psi,\phi)-\sigma(\phi,\psi)> + Re<\phi,\partial_A\phi +\alpha \psi>. 
\end{equation}
The first term is the $L^2$-inner product of forms and the second is the
$L^2$-inner product of sections of $S\otimes L$.

The quadratic form (\ref{2forma}) induces an operator on the
$L^2_k$-tangent bundle of ${\cal B}_H$, which is the same as the
linearization $T$ of the flow equations. At a critical point this is
the Hessian of the functional $\tilde {\cal C}$.
In fact the tangent space of ${\cal B}$ at $(A,\psi)$ is given by the 
pairs $(\alpha,\phi)$ that satisfy $G^*(\alpha,\phi)=0$.
In fact, the following holds.

\begin{lem}
the following relation is satisfied for all $(\alpha,\phi)$:
\[ G^*_{A,\psi}( T_{(A,\psi)}(\alpha,\phi))=0. \]
\end{lem}

\noindent\underline{Proof:} We have
\[ G^*_{A,\psi}( T_{(A,\psi)}(\alpha,\phi))= -d^*(*d\alpha
-\sigma(\psi,\phi) -\sigma(\phi,\psi)) + \]
\[ + i Im <\psi,\partial_A\phi +\alpha\psi >=d^*(\sigma(\psi,\phi)
+\sigma(\phi,\psi))+i Im <\psi,\partial_A\phi>= \]
\[ = -\frac{1}{2}<\psi,\partial_A\phi> +\frac{1}{2}< \partial_A\phi,
\psi> +i Im <\psi,\partial_A\phi>=0, \]
using the fact that $\partial_A\psi=0$ and that 
\[ d^*\sigma(\psi,\psi)=
\frac{1}{2} (<\partial_A\psi,\psi>-<\psi,\partial_A\psi>). \]

\noindent QED

Since $(A,\psi)$ is a point in ${\cal M}_C$, by elliptic regularity we
can regard the operator $T$ at a critical point as an operator that
maps the $L^2_k$-tangent space to itself. This completes the proof of
Theorem \ref{hessian}.

\noindent QED

\section{Homology}

In this section we shall assign a {\em Morse index} to
the critical points of $\tilde {\cal C}$, upon fixing the Morse index
of one particular critical point.  
In order to construct the Floer homology, we shall consider paths of
steepest descent that connect critical points of relative Morse index
equal to one. The boundary operator of the Floer homology is
constructed by counting these paths with their orientation as in \cite{Fl}.

\subsection{Gradient Flow}

The analysis presented in this and the following two sections is part
of a joint work with B.L. Wang and will appear in \cite{MaWa}.

We introduce suitable moduli spaces of gradient flow lines connecting
critical points. We prove in the following that, generically, these
are smooth manifolds that are cut out transversely, hence with the
dimension prescribed by the index theorem. This property depends on an
accurate choice of a class of perturbations for the gradient flow
equations. 

Consider the space connections and sections $(\AAA,\Psi)$ on $Y\times
\RR$ with the product metric $g+dt^2$, topologized with the
weighted Sobolev norms \cite{LM} \cite{Fl}. Here we choose the weight
$e_\delta(t)=e^{\tilde\delta t}$, where $\tilde\delta$ is a smooth
function with bounded derivatives, $\tilde\delta : \RR\to
[-\delta,\delta]$ for some fixed positive number $\delta$, such that
$\tilde\delta (t)\equiv -\delta$ for $t\leq -1$ and $\tilde\delta
(t)\equiv \delta$ for $t\geq 1$. 
The $L^2_{k,\delta}$ norm is defined as $\| f \|_{2,k,\delta}=\|
e_\delta f \|_{2,k}$; the weight $e_\delta$ imposes an exponential
decay as asymptotic condition along the cylinder.

\begin{prop}
Let $Y$ be a compact oriented three-manifold endowed with a fixed
Riemannian metric $g_0$. Consider the cylinder $Y\times \RR$ with the
metric $g_0+dt^2$.
The weighted Sobolev spaces $L^2_{k,\delta}$ on the manifold $Y\times
\RR$ satisfy the following Sobolev embeddings.

(i) The embedding $L^2_{k,\delta}\hookrightarrow L^2_{k-1,\delta}$ is
compact for all $k\geq 1$.

(ii) If $k> m+2$ we have a continuous embedding
$L^2_{k,\delta}\hookrightarrow {\cal C}^m$.

(iii) If $k> m+3$ the embedding $L^2_{k,\delta}\hookrightarrow {\cal
  C}^m$ is compact.

(iv) If $2<k'$ and $k\leq k'$ the multiplication map 
$L^2_{k,\delta}\otimes L^2_{k,\delta}\stackrel{m}{\to} L^2_{k,\delta
  /2}$ is continuous.

\noindent Consider a metric $g_t +dt^2$ on the cylinder $Y\times
\RR$ such that for a fixed $T$ we have $g_t\equiv g_0$ for $t\geq T$
and $g_t\equiv g_1$ for $t\leq -T$ and $g_t$ varies smoothly when
$t\in [-1,1]$. The same Sobolev embedding theorems hold for the
$L^2_{k,\delta}$ spaces on $(Y\times \RR,g_t +dt^2)$.
\label{sobolev}
\end{prop}

Choose smooth representatives $(A_0,\psi_0)$ and $(A_1,\psi_1)$ of $a$
and $b$ in ${\cal M}$. Choose a smooth path $(A(t),\psi(t))$ such that
for $t\leq 0$ it satisfies $(A(t),\psi(t))\equiv (A_0,\psi_0)$ and for
$t\geq 1$ it is $(A(t),\psi(t))\equiv (A_1,\psi_1)$.

Let ${\cal A}_{k,\delta}(a,b)$ be the space of pairs $(\AAA,\Psi)$ on
$Y\times \RR$ satisfying
\[ (\AAA,\Psi)\in (A(t),\psi(t))+L^2_{k,\delta} (\Lambda^1(Y\times\RR)\oplus
\Gamma (S^+\otimes L)). \]

Consider the group ${\cal G}_{k+1,\delta}$ of gauge transformations in
${\cal G}(Y\times \RR)$, locally modelled on
$L^2_{k+1,\delta}(\Lambda^0(Y\times\RR)$, that decay to asymptotic
values $\lambda_{\pm\infty}$ which satisfy condition
(\ref{int=0}) on $Y$. This gauge group acts on ${\cal
  A}_{k,\delta}(a,b)$ and we can form the quotient ${\cal
  B}^H_{k,\delta}(a,b)$.

We consider the perturbed gradient flow equations for a path
$(A(t),\psi(t))$,
\beq
\frac{d}{dt}\psi(t)=-\partial_{A(t)}\psi (t)
\label{3SW1P'}
\eeq
and
\beq
\frac{d}{dt}A(t)=\sigma(\psi(t),\psi(t))-*F_{A(t)}+2i\rho
+2q_{(\AAA,\Psi)}(t).
\label{3SW2P'}
\eeq

Equations (\ref{3SW1P'}) and (\ref{3SW2P'}) can be rewritten in terms
of pairs $(\AAA,\Psi)$ in the form
\beq
D_{\AAA}\Psi=0
\label{4SW1P'}
\eeq
and
\beq
F_{\AAA}^+=\Psi\cdot\bar\Psi +i\mu +P_{(\AAA,\Psi)}.
\label{4SW2P'}
\eeq

The perturbation $P=*q+q\wedge dt$ is a function of ${\cal B}^H(a,b)$
to $\pi^*(\Lambda^1(Y))\cong\Lambda^{2+}(Y\times\RR)$, such that the
corresponding equations in a temporal gauge
(\ref{3SW1P'}) and (\ref{3SW2P'}) are preserved under the action of $\RR$
by reparametrizations of the path $(A(t),\psi(t))$. The class of such
perturbations is described as follows.

\begin{defin}
\label{calP}
The space of perturbations ${\cal P}$ is the space of maps
\[ P:{\cal B}^H_{k,\delta}(a,b)\to
L^2_{k,\delta}\Lambda^{2+}(Y\times\RR), \] 
that satisfy the following conditions.

\noindent (1) $P_{(\AAA,\Psi)}= *q_{(\AAA,\Psi)}(t)
+q_{(\AAA,\Psi)}(t) \wedge dt$, where $q_{(\AAA,\Psi)}(t)$
satisfies 
\[ q_{(\AAA,\Psi)^T}(t)=q_{(\AAA,\Psi)}(t+T), \]
where $(\AAA,\Psi)^T$ is the $T$-translate of $(\AAA,\Psi)$, namely the
pair that is represented in a temporal gauge by $(A(t+T),\psi(t+T))$.

\noindent (2) the $L^2_{k,\delta}$-norm of the perturbation
$P_{(\AAA,\Psi)}$ is bounded uniformly with respect to $(\AAA,\Psi)$;

\noindent (3) the linearization ${\cal D}P_{(\AAA,\Psi)}$ is a compact
operator from the $L^2_{k,\delta}$ to the
$L^2_{k-1,\delta}$ tangent spaces. 

\noindent (4) for all $l\leq k-1$, we have that 
\[ \| q_{(\AAA,\Psi)}(t)\|_{L^2_l}\leq C_l \| \nabla \tilde{\cal
    C} (A(t),\psi(t)) \|_{L^2_l} \]
in the $L^2_l$-norm on $Y\times \{ t\}$, for all $|t|>T$, where
\beq
\nabla \tilde{\cal C} (A(t),\psi(t))=(-\partial_{A(t)}\psi (t),
\sigma(\psi(t),\psi(t))-*F_{A(t)}+2\rho)
\eeq
is the gradient flow of the functional $\tilde{\cal C}$, and $0<C_l<1$;

\end{defin}

With a perturbation in the class ${\cal P}$ the equations
(\ref{3SW1P'}) and (\ref{3SW2P'}) are invariant 
with respect to the action of $\RR$ by translations along the gradient
flow lines, that is if $(A(t), \psi(t))$ is a solution of (\ref{3SW1P'}) and
(\ref{3SW2P'}), then $(A(t+T), \psi(t+T))$ is also a solution for any
$T\in \RR$. 

An example of 
perturbation with these properties has been constructed by
Froyshov \cite{Fr}. 

\begin{prop}
The class of perturbations introduced by Froyshov in \cite{Fr} is in our
class ${\cal P}$. 
\label{froyshov}
\end{prop}

According to Froyshov's construction, for fixed smooth compactly
supported functions $\eta_1$, $\eta_2$, with 
$supp(\eta_1)\subset [-1,1]$ and
$\eta_2|_I(t)=t$ on an interval $I$ containing all the critical values
of $\tilde{\cal C}$, a function $h: {\cal B}^H_{k,\delta}(a,b)\to {\cal
  C}^m (\RR)$ is defined as 
\[ h_{(\AAA,\Psi)}(T)= \int_{\RR}\eta_1(s-T) \eta_2
(\int_{\RR}\eta_1(t-s) \tilde{\cal C}(A(t),\psi(t))dt)ds. \]

Let $\Lambda_\Xi^2(Y\times\RR)$ be the set of ${\cal C}^m$ 2-forms $\omega$ 
that are compactly supported in $Y\times\Xi$, where $\Xi$ is the
complement of a union of small intervals centered at the critical
values of $\tilde{\cal C}$. Froyshov's perturbation is constructed by setting
\[ P_{(\AAA,\Psi)}=(h_{(\AAA,\Psi)}^*(\omega))^+, \]
where $h_{(\AAA,\Psi)}^*(\omega)$ is the pullback of $\omega$ along the
map $Id_Y\times h_{(\AAA,\Psi)}:Y\times \RR\to Y\times \RR$.

As shown in \cite{Fr}, the function $h_{(\AAA,\Psi)}$ is bounded with
all derivatives, uniformly with respect to 
$(\AAA,\Psi)$. Moreover, by the choice of $\Xi$, the perturbation
$h_{(\AAA,\Psi)}^*(\omega)$ is smooth and compactly supported, hence in
$L^2_{k,\delta}$.
  
Condition (1) holds, since the function $h_{(\AAA,\Psi)}(t)$ satisfies 
\[ h_{(\AAA,\Psi)}(t+\tau)=h_{(\AAA,\Psi)^\tau}(t), \]
where $(\AAA,\Psi)^\tau$ is the $\tau$-reparametrized solution
represented in a temporal gauge by $(A(t+\tau),\psi(t+\tau))$.
In fact,
\[ h_{(\AAA,\Psi)^\tau}(T)=\int_{\RR}\eta_1(s-T) \eta_2
(\int_{\RR}\eta_1(t-s) \tilde{\cal C}(A(t+\tau),\psi(t+\tau))dt)ds= \] 
\[ =\int_{\RR}\eta_1(s-T) \eta_2(\int_{\RR}\eta_1(u-s-\tau) \tilde{\cal
  C}(A(u),\psi(u))du)ds= \]
\[ = \int_{\RR}\eta_1(v-T-\tau) \eta_2(\int_{\RR}\eta_1(u-v) \tilde{\cal
  C}(A(u),\psi(u))du)dv=h_{(\AAA,\Psi)}(T+\tau). \]

Condition (2) holds: in fact, it is shown in \cite{Fr} that the function
$h_{(\AAA,\Psi)}$ is bounded with all derivatives, uniformly with respect to 
$(\AAA,\Psi)$. The Sobolev embeddings of Proposition \ref{sobolev}
provide the uniform bound in the $L^2_{k,\delta}$-norms.

Condition (3) also follows from Froyshov (\cite{Fr}, Prop.5): for
$\omega$ a ${\cal C}^m$ form, the linearization of the perturbation
$h_{(\AAA,\Psi)}^*(\omega)$ at the point $(\omega,\AAA,\Psi)$ is a
  bounded operator $K_{(\omega,\AAA,\Psi)}:L^2_{k,\delta}\to {\cal C}^m$ with 
\[ supp\left(K_{(\omega,\AAA,\Psi)}(\alpha,\Phi)\right)\subset
h^{-1}_{(\AAA,\Psi)}(\Xi)\times Y. \]

Condition (4) follows from the fact that the perturbation
$h_{(\AAA,\Psi)}^*(\omega)$  is compactly supported. That is,
$q_{(\AAA,\Psi)}(t)$ will be identically zero for large enough $T$ and
in particular (4) is satisfied for large enough $t$. 

Notice that in properties (3) and (4), the interval $[-T,T]$ cannot be
chosen uniformly with respect of $(\AAA,\Psi)$.
 
Other suitable 
perturbations of the functional ${\cal C}$ can be used to achieve
transversality of the moduli space of flow lines. For instance one can
consider perturbations by functions of the holonomy of the connection
$A$. This kind of perturbation is used in \cite{CMW}.

\subsection{Flow lines and transversality}

Let ${\cal L}_{(\AAA,\Psi)}$ be the linearization of equations
(\ref{4SW1P'}) and (\ref{4SW2P'}) on ${\cal B}^H_{k,\delta}(a,b)$.

The operator ${\cal L}$ is of the form 
\[ {\cal L}_{(\AAA,\Psi,P)}(\alpha,\Phi)=\left\{\begin{array}{c}
D_{\AAA}\Phi+ \alpha \Psi \\ d^+\alpha
-\frac{1}{2}Im(\Psi\cdot\bar\Phi)+{\cal D}P_{(\AAA,\Psi)}(\alpha,\Phi)\\
G^*_{(\AAA,\Psi)}(\alpha,\Phi)
\end{array}\right. \]
mapping
\[ L^2_{k,\delta}(\Lambda^1(Y\times\RR)\oplus \Gamma(S^+\otimes
L))\to L^2_{k-1,\delta} (\Lambda^0(Y\times\RR)\oplus \Lambda^{2+}
(Y\times\RR). \]
The operator $G^*$ is the adjoint of 
$G_{(\AAA,\Psi)}(f)=(-df,f\Psi)$.

As the following proposition shows,
the operator ${\cal L}_{(\AAA,\Psi,P)}$ is obtained by adding the small
perturbation ${\cal D}P_{(\AAA,\Psi)}$ to a Fredholm map
from $L^2_{k,\delta}$ to $L^2_{k-1,\delta}$, hence it is
Fredholm. Therefore 
we have a well defined relative Morse index of two critical points $a$
and $b$ in ${\cal M}^H_C$.

\begin{prop}
Suppose $a$ and $b$ are irreducible critical points for the functional
$\tilde{\cal C}$. Let $\{ \lambda_a \}$ and $\{ \lambda_b \}$ be the
eigenvalues of the Hessian $T$ at the points $a$ and $b$.
Assume that the positive number $\delta$ satisfies 
$\delta < \min \{ |\lambda_a|, |\lambda_b| \}$. 
Let $(\AAA,\Psi)$ be a solution of (\ref{4SW1P'}) and (\ref{4SW2P'}) in
${\cal B}_{k,\delta}(a,b)$. Then the linearization ${\cal
  L}_{(\AAA,\Psi)}$ is a Fredholm operator of index
\[ Ind({\cal L}_{(\AAA,\Psi)})=\sigma(a,b). \]
The right hand side $\sigma(a,b)$ is the spectral flow of the operator
$\nabla{\cal F}$ along a path $(A(t),\psi(t))$ in ${\cal A}$
that corresponds to $(\AAA,\Psi)$ under $\pi^*$. 
The quantity $\sigma(a,b)$ is independent of the path, hence
\[ \sigma(a,b)=\mu(a)-\mu(b) \]
defines a relative Morse index of $a$ and $b$, where $\mu(a)$ is the
spectral flow of $\nabla{\cal F}$ on a path joining $a$ to a fixed
$[A_0,\psi_0]$ in ${\cal M}^H_C$.
\label{relmorseind}
\end{prop}

The Fredholm property follows from \cite{LM} theorem 1.3.
The result about the spectral flow requires the following lemma proven
by R.G. Wang \cite{GWa}. 

\begin{lem}
Suppose given a path $(A(t),\psi(t))$ in ${\cal A}$ that decays exponentially
fast to asymptotic values $(A,\psi)$ and $\lambda (A,\psi)$ in the
same gauge class, with $\lambda\in {\cal G}$. Then the index of the
linearization 
\[ Ind({\cal L}_{(\AAA,\Psi)})=\frac{i}{2\pi}\int_Y c_1(L)\wedge
\lambda^{-1}d\lambda. \]
\label{spectralsum}
\end{lem}

We report here the proof given in \cite{GWa}. Namely, $\lambda : Y\to
U(1)$ determines a $U(1)$ bundle over $Y\times S^1$, by identifying
the ends of the cylinder; the connection $A(t)$ gives rise to a
connection on this line bundle $\hat L$ over $Y\times S^1$. The index
of the linearization ${\cal L}_{(\AAA,\Psi)}$ is therefore given by
\[ Ind({\cal L}_{(\AAA,\Psi)})=-\frac{1}{16\pi^2}\int_{Y\times S^1}
c_1(\hat L)^2-\frac{2\chi+3\sigma}{4}, \]
since ${\cal L}_{(\AAA,\Psi)}$ is, up to compact perturbations, the
linearization of the four-dimensional Seiberg-Witten equations on
$Y\times S^1$ with the $Spin_c$ structure specified by the line bundle
$\hat L^2$. The term $2\chi+3\sigma=0$ on a manifold of the form
$Y\times S^1$. As for the first term, notice that we can write 
\[ F_{A(t)}\wedge F_{A(t)}=F_{A(t)}\wedge \frac{dA(t)}{dt} \wedge dt, \] 
and therefore we get
\[ \frac{-1}{8\pi^2}\int_{Y\times S^1}F_{A(t)}\wedge \frac{dA(t)}{dt}
\wedge dt= \frac{i}{2\pi} \int_Y c_1(L)\int_{S^1} \frac{dA}{dt}= \]
\[ = \frac{i}{2\pi} \int_Y c_1(L)\wedge \lambda^{-1}d\lambda, \]
since $A(+\infty)-A(-\infty)=\lambda^{-1}d\lambda$.

The spectral flow of the family of operators
$T_t=T_{(A(t),\psi(t))}$, which is the index of the operator 
$\frac{\partial}{\partial t}+ T_{(A(t),\psi(t))}$ by \cite{APS}, p.57,
p.95. Suppose given a path $(A(t),\psi(t))$ in ${\cal A}$ with
endpoints that are gauge equivalent through a gauge transformation in
$\tilde{\cal G}$. Then, by the previous Lemma, the spectral flow is given by
\[ Ind(\frac{\partial}{\partial t}+ T_{(A(t),\psi(t))})=0, \]
that is, the spectral flow around a loop in ${\cal B}^H$ is trivial.

The additivity of spectral flows implies that 
\[ \sigma(a,c)=\mu(a)-\mu(c)=\sigma(a,b)+\sigma(b,c). \]

This means that the relative Morse index of critical points in ${\cal
  M}^H_C$ is well defined and equal to the spectral flow.
We can consider the same construction of moduli spaces of flow lines
in ${\cal B}_{k,\delta}(a,b)$, that is modulo the action of the full
gauge group ${\cal G}(Y\times \RR)$. This determines a relative Morse
index for points in ${\cal M}_C$ given by the spactral flow. However,
in this case the relative Morse index is only defined up to an integer
multiple of $l$, where
\beq
l=g.c.d.\{ <c_1(L)\cup h,[Y]> ~| h\in H \}.
\label{ell}
\eeq
This follows from the spectral formula \ref{spectralsum}. Notice that
$l$ is an even number.

Thus we have an important difference between 
Seiberg--Witten and Donaldson Floer homology: 
the relative index is well defined and
there is no ambiguity coming from loops in $\hat{\cal B}^0$. Thus in
our case the Floer groups will be $\ZZ$-graded. This makes the
Seiberg-Witten-Floer homology more similar to the finite dimensional
case \cite{W1}. In fact, as proved in \cite{MST}, for
three-manifolds of the form $\Sigma\times S^1$, with $\Sigma$ a
Riemann surface, the Floer groups are just the ordinary homology of
a symmetric product of $\Sigma$. Similar results for the mapping cylinder
of a surface $\Sigma$ have been obtained by A.L. Carey and B.L. Wang 
\cite{CaWa}. 

The ambiguity mod 8 in Donaldson Floer homology is related 
to the possibility of rescaling and gluing instantons at different
$t\in \RR$ on $Y\times \RR$, as explained in \cite{Fl}. There is a
form of periodicity in the Seiberg-Witten-Floer groups as well if we
consider the moduli space ${\cal M}_C$ instead of ${\cal M}_C^H$, as
we discuss later. However, this depends only upon the nature of the
covering ${\cal B}^H\to {\cal B}$. In fact, the non-rescaling property
of the solutions of the monopole equations imply the absence of
``bubbling'' phenomena. 

Consider the moduli space ${\cal M}^H(a,b)$ of solutions of the
equations (\ref{4SW1P'}) and (\ref{4SW2P'}) in ${\cal
  B}^H_{k,\delta}(a,b)$. 

\begin{prop}
Given $a$ and $b$, two critical points of $\tilde{\cal C}$,
for a generic choice of the perturbation $P\in {\cal P}$, the moduli space 
${\cal M}^H(a,b)$ of gradient flow lines is a smooth manifold, cut out
transversely by the equations, of dimension 
\[ \dim({\cal M}^H(a,b))=\mu(a)-\mu(b), \]
where $\mu(a)-\mu(b)$ is the relative Morse index of the critical points.
\label{transverse}
\end{prop}

\noindent\underline{Proof:}
It is first necessary to prove that there are no reducible flow lines
connecting the critical points $a$ and $b$. This is proven later in
the next section. Once it is known that ${\cal M}^H(a,b)$ lies in the
irreducible component $\hat{\cal B}^H_{k,\delta}(a,b)$, the statement
follows via the implicit function theorem, upon showing that, for a
generic choice of the perturbation $P$, the linearization ${\cal L}$
is surjective.   

Consider the operator
\[ \hat{\cal L}_{(\AAA,\Psi,P)}(\alpha,\Phi,p)={\cal
  L}_{(\AAA,\Psi,P)}(\alpha,\Phi) +p_{(\AAA,\Psi,P)}(\alpha,\Phi), \]
where we vary the perturbation by an element $p_{(\AAA,\Psi,P)}$ of
$T_P{\cal P}$. This corresponds to varying the parameter
$\omega\in\Lambda_\Xi^2(Y\times \RR)$ in Froyshov's class of
perturbations. 

The operator ${\cal L}$ is Fredholm, therefore $\hat{\cal L}$ has a
closed range. We show that $\hat{\cal L}$ is surjective by proving
that it has dense range.
 
Suppose given an element $(\beta,\xi,g)$ in
$L^2_{-k-1,-\delta}$ that is $L^2$-orthogonal 
to the range of the operator $\hat{\cal L}$. Here $\beta$ is an element in 
$\Lambda^{2+}(Y\times\RR, i\RR))$, $\xi$ is a spinor, and $g$ is a
zero-form. The element $(\beta,\xi,g)$ is in the kernel of the adjoint
$\tilde{\cal L}^*$, which is an elliptic operator with
$L^2_{-k,-\delta}$ coefficients, thus $(\beta,\xi)$ lives in
$L^2_{-k,-\delta}$ by elliptic regularity. 
For the same reason, $g$ lives in $L^2_{k-1,\delta}$, since the
operator $G^*$ is the adjoint of $G$ with respect to the
$L^2_{k-1,\delta}$ inner product. 
If we consider the $L^2$-pairing of $L^2_{k,\delta}$ and $L^2_{-k,-\delta}$,
we get 
\[ \langle \beta, d^+\alpha -\frac{1}{2}Im(\Psi\cdot\bar\Phi))+{\cal
  D}P_{(\AAA,\Psi)}(\alpha,\Phi)+p_{(\AAA,\Psi,P)}(\alpha,\Phi) \rangle
+ \] 
\[ +\langle \xi, D_{\AAA}\Phi+\alpha\Psi \rangle+ 
\langle g,G^*_{(\AAA,\Psi)}(\alpha,\Phi)\rangle=0. \]

By varying $p\in {\cal P}$ we force $\beta\equiv 0$.
The remaining inner product 
\[ \langle \xi, D_{\AAA}\Phi+\alpha\Psi \rangle+
\langle g,G^*_{(\AAA,\Psi)}(\alpha,\Phi)\rangle=0 \]
gives the following equations

(a) $(e_{-\delta}de_{\delta})g=\frac{1}{2}\xi\cdot\bar\Psi$ and

(b) $D_A\xi-g\Psi=0$.

\noindent We assume that $\Psi$ is not identically zero. Applying 
$d^*$ to (a) and using (b) we obtain
$d^*(e_{-\delta}de_{\delta}g)+g|\Psi|^2=0$. Equivalently, we get
\[ (e_{\delta/2}d^*e_{-\delta/2}) (e_{-\delta/2}de_{\delta/2})
e_{\delta/2}g +|\Psi|^2 e_{\delta/2}g=0. \] 
The equation 
\[ \Delta_{\delta/2} g + e_{\delta/2}g|\Psi|^2=0, \] 
with
\[ \Delta_{\delta}=e_{-\delta}\Delta e_{\delta}, \] 
omplies that $g\equiv 0$, since $g$ decays at $\pm\infty$ and the
maximum priciple applies. Then, by varying 
$\alpha$ alone in $\langle \xi, D_{\AAA}\Phi+\alpha\Psi \rangle=0$,
we force $\xi$ to vanish on some open
set. The pair $(\beta,\xi)$ is in the kernel of the operator
$T_{(\AAA,\Psi)}^*$, hence by analytic continuation $(\beta,\xi)\equiv 0$.

Thus the operator $\tilde{\cal L}$ is surjective. This implies that
zero is a regular value for the map defined by the equations
(\ref{4SW1P'}) and (\ref{4SW2P'}). Therefore the moduli space  ${\cal
  M}od$ of triples $([\AAA,\Psi],P)$ in $\hat{\cal B}^H_{L^2_{k,\delta}}(a,b)
  \oplus {\cal P}$ that satisfy the equations is a 
smooth (infinite dimensional) manifold with virtual tangent space
$Ker(\hat{\cal L})$. 

The projection $\Pi: {\cal M}od \to {\cal P}$ given by
$\Pi([\AAA,\Psi],P)=P$ linearizes to a surjective Fredholm operator
$${\cal D}\Pi:Ker(\tilde{\cal L})\to T_P{\cal P}.$$ 
The kernel of ${\cal
  D}\Pi$ is $Ker({\cal D}\Pi_{(\AAA,\Psi,P)})=Ker({\cal L}_{(\AAA,\Psi,P)})$.
The infinite dimensional Sard
theorem implies that the moduli space ${\cal M}^H(a,b)$, 
for a generic perturbation $P\in {\cal P}$, is
the inverse image under the projection map from ${\cal M}od$ to ${\cal
  P}$ of a regular value. Thus ${\cal M}^H(a,b)$ is a
smooth manifold which is cut out transversely by the equations.
Equivalently, the linearization ${\cal L}$ with a fixed generic $q$ is
surjective.  

As we shall prove in the following, the virtual dimension of the
moduli space ${\cal M}^H(a,b)$ equals the  
index of the Fredholm operator ${\cal L}$. According to Proposition
\ref{relmorseind}, this is the relative Morse index $\mu(a)-\mu(b)$.

\noindent QED

The case of points $a$ and $b$ in ${\cal M}_C$ is analogous. In this
case we obtain a moduli space ${\cal M}(a,b)$ with components of
dimension $\mu(a)-\mu(b)+\xi(h)$ with $h\in H$, according to the
spectral formula \ref{spectralsum}.

\subsection{Flow equations}

It is clear that moduli spaces of flow lines can be defined either
using equations (\ref{3SW1P'}) and (\ref{3SW2P'}), for pairs
$(A(t),\psi(t))$ in a temporal gauge with gauge tranformations in
${\cal G}$, or using the corresponding equations (\ref{4SW1P'}) and
(\ref{4SW2P'}) for pairs $(\AAA,\Psi)$ with gauge action of ${\cal
  G}(a,b)$, as we did in the previous section.

The two descriptions are equivalent,
however, the Fredholm analysis is somewhat simpler in the latter
case. If fact, in order to set up a good Freholm analysis with
equations (\ref{3SW1P'}) and (\ref{3SW2P'}), one has to add a
correction term $\gamma(A,\psi)$ in the flow equations in order to
make the flow tangent to a fixed slice of the gauge action at a point
$[A_0,\psi_0]$. This corresponds to the analysis worked out in
\cite{MMR} for the case of Donaldson Floer theory. One difficulty
arises in this case, since the correction term
$\gamma(A,\psi)$ does not preserve the temporal gauge condition. This
problem can be overcome by replacing the temporal gauge condition with
the condition of {\em standard form} introduced in \cite{MMR} and
allowing time-dependent gauge transformations.

The linearization of the equations contains the extra
term $L_\gamma$ that linearizes $\gamma(A,\psi)$,
\[ {\cal L}_{(A(t),\psi(t),q)}(\alpha,\phi)=\frac{\partial}{\partial t}+
L_{A(t),\psi(t)} +{\cal D}q_{(A(t),\psi(t))}+L_\gamma. \]  
Anisotropic Sobolev norms $L^2_{(k,m)}$ can be introduced on the spaces
of connections and sections and gauge transformations, as analysed in
\cite{MMR}.The linearization $L_\gamma$ is a compact operator
with respect to these norms. We have the equivalent of Proposition
\ref{relmorseind} within this formulation: ${\cal
  L}_{(A(t),\psi(t),q)}$ is a Fredholm operator whose index is the
spectral flow of the family of operators $L_{(A(t),\psi(t))}$. 
The transversality result can also be proven in this context. This
formulation, however, will not be worked out in this paper.

The key to the equivalence of these two formulations is the elliptic
regularity and a decay estimate that will be proven later.
In fact, if we choose a smooth perturbation in ${\cal
  P}$, by elliptic regularity it is possible to represent any
solution of (\ref{4SW1P'}) and (\ref{4SW2P'}) in ${\cal
  B}_{k,\delta}(a,b)$ by a smooth representative. This defines a
solution of (\ref{3SW1P'}) and (\ref{3SW2P'}) in $L^2_{(k,m)}$ in a
standard form. Conversely solutions of (\ref{3SW1P'}) and
(\ref{3SW2P'}) give rise to solutions of (\ref{4SW1P'}) and
(\ref{4SW2P'}) in ${\cal B}_{k,\delta}(a,b)$ due to the exponential
decay towards the endpoints that will be proven in the next section.

\subsection{Decay estimate}

In this subsection we show that the moduli space ${\cal M}^H(a, b)$
only contains irreducible flow lines. Given that $a$ or $b$ are
irreducible points, the result follows directly from the decay estimate
below. The analysis is based on \cite{MST} and \cite{Wa}.

\begin{thm}
Let $a$ and $b$ be non-degenerate critical points in ${\cal
  B}^H$. There exists a weight $\delta >0$ such that the following holds. 
Suppose given any solution
$[\AAA,\Psi]$ of (\ref{4SW1P'}) and (\ref{4SW2P'}) that is represented
by a smooth pair $(A(t),\psi(t))$ in a temporal gauge, with asymptotic
values $(A_a,\psi_a)$ and $(A_b,\psi_b)$ representing the elements $a$
and $b$. Then there exists a constant $K$ such that, for $t$
outside an interval $[-T,T]$, the distance in any fixed
${\cal C}^l$-topology of $(A(t),\psi(t))$ from the
endpoints is
\[ dist_{{\cal C}^l}((A(t),\psi(t)),(A_i,\psi_i))< K \exp(-\delta
|t|), \]
with $i=a$ if $t< -T$ and $i=b$ if $t>T$. 
\label{decay}
\end{thm}

\noindent\underline{Proof:}
The proof consists of a few steps, mainly based on the analysis worked
out in \cite{MST}. Let us consider the decaying as $t\to \infty$; the
other case is analogous. 

{\bf Claim 1:}
Let $(A(t),\psi(t))$ be a solution of the flow equations
(\ref{3SW1P'}) and (\ref{3SW2P'}) with finite energy and with
limit $\lim_{t\to\infty}(A(t),\psi(t))=(A_b,\psi_b)$. Then there is a
$T>>0$ and a constant $K_b$, such that
the inequality
\[  \tilde{\cal C}(A(t),\psi(t))-\tilde{\cal C}
(A_b,\psi_b) \leq K_b  
\| \nabla \tilde{\cal C} (A(t),\psi(t)) \|_{L^2}^2 \]
holds for $t\geq T$. 

Lojasiewicz inequality \cite{Si} shows that there exists $T>0$ and an
exponent $0<\theta<1/2$ such that
\[ |\tilde{\cal C} (A(t),\psi(t))-\tilde{\cal C} ((A_b,\psi_b))
|^{1-\theta}\leq K_b \|\nabla \tilde{\cal C}(A(t),\psi(t)) \|_{L^2}, \]
for $t\geq T$. 
Under the assumption of non-degenerate Hessian at the point $b$, the
exponent can be improved to $\theta=\frac{1}{2}$.

{\bf Claim 2:}
For a solution $(A(t),\psi(t))$ of (\ref{3SW1P'}) and (\ref{3SW2P'}), 
the inequality 
\[ \frac{1}{2} \int_t^\infty \|\nabla \tilde{\cal C}(A(s),\psi(s)) \|^2 ds
\leq \tilde{\cal C} (A(t),\psi(t))- \]
\[ -\tilde{\cal C}(A_b,\psi_b) \leq \frac{3}{2} \int_t^\infty \|\nabla
\tilde{\cal C}(A(s),\psi(s)) \|^2 ds \]
holds for large $t$. 

Without loss of generality we can assume that the perturbation in
${\cal P}$ satisfies Condition (4) of \ref{calP} with $C_0< 1/2$,
so that
\[ \|q_{(\AAA,\Psi)}(t)\|_{L^2}< \frac{1}{2}\| \nabla\tilde{\cal C}
(A(t),\psi(t)) \|_{L^2}. \]
Thus, we can replace the equality
\[ \tilde{\cal C} (A(t),\psi(t))-\tilde{\cal C}
((A_b,\psi_b))= 
\int_t^\infty -\frac{d}{ds}\tilde{\cal C} (A(s),\psi(s))
ds \]
\[ =\int_t^\infty -<\frac{d}{ds} (A(s),\psi(s)),\nabla \tilde{\cal C}
(A(s),\psi(s))> ds = \]
\[ =\int_t^\infty \|\nabla \tilde{\cal C} (A(s),\psi(s)) \|^2 ds, \] 
that holds for solutions of the unperturbed equations with the
inequality of Claim 2 for solutions of the perturbed equations.

{\bf Claim 3:} The quantity 
\[ E(t)=\int_t^\infty \|\nabla\tilde{\cal C} (A(s),\psi(s)) \|^2 ds \]
decays exponentially as $t\to \infty$.

In fact, the inequality of Claim 2 gives the first inequality in the following
estimate: 
\[ E(t)\leq 2 (\tilde{\cal C} (A(t),\psi(t))-\tilde{\cal
  C}(A_b,\psi_b)) \leq \]
\[  \leq K_b \| \nabla \tilde{\cal C} (A(t),\psi(t))
\|^2  = -K_b \frac{d}{dt} E(t). \]
The second inequality follows fron Claim 1 with the best exponent
$\theta=1/2$.  

{\bf Claim 4:} For large $t$ we have the inequality 
\[ dist_{L^2_1}((A(t),\psi(t)),(A_b,\psi_b))\leq K \int_{t-1}^\infty
\| \nabla \tilde{\cal C} (A(s),\psi(s)) \|_{L^2}^2 ds. \]

In fact, the perturbation in ${\cal P}$ satisfies
\[ \| q_{(A(t),\psi(t))}(t)\|_{L^2_1} \leq C_1 \|
\nabla \tilde{\cal C}(A(t),\psi(t)) \|_{L^2_1}, \]
with $0< C_1 <1$ for large $t$. Thus we have
\[ dist_{L^2_1}((A(t),\psi(t)),(A_b,\psi_b))\leq K \int_{t}^\infty
\| \nabla \tilde{\cal C} (A(s),\psi(s)) \|_{L^2_1}^2 ds. \]
Lemma 6.14 of \cite{MST} implies that the latter is bounded by
\[ K\int_{t-1}^\infty \| \nabla \tilde{\cal C} (A(s),\psi(s))
\|_{L^2}^2 ds, \] thus proving the inequality. 

The exponential decay of $E(t)$ proves the claim of the theorem for
the case of $L^2_1$-topology. Smooth estimates then follow by a
bootstrapping argument and elliptic regularity. 

\noindent QED

Analogous exponential decay estimates have been proven in
\cite{Wa}. An immediate corollary of Theorem \ref{decay} is the
following. 

\begin{corol}
No reducible solution arises among the flow lines in ${\cal M}^H(a,b)$. 
\end{corol}

In fact, $(A_a,\psi_a)$ and $(A_b,\psi_b)$ have non-trivial
spinor. Thus, by the exponential decay estimate, $\psi(t)$ is forced
to be non-trivial.

\subsection{The Boundary Operator}

We have described all the ingredients that are needed in order to 
construct the Floer homology. Following \cite{Fl} we can define
the boundary operator.

Let $\hat{\cal M}^H(a,b)$ and $\hat{\cal M}(a,b)$ be the quotients of
${\cal M}^H(a,b)$ and of ${\cal M}(a,b)$ by the action of $\RR$ by
translations. 

\begin{defin}
Let $a$ and $b$ be two critical points in ${\cal M}_C^H$ of relative 
Morse index $\mu(a)-\mu(b)=1$. Define the boundary $\partial$ to be
the operator with matrix elements
\[ <\partial^H a, b>=\epsilon^H(a,b), \]
where $\epsilon^H(a,b)$ is the algebraic sum over the paths joining $a$
and $b$ of the signs given by the orientation,
\[ \epsilon^H(a,b)=\sum_{\gamma\in\hat{\cal M}^H(a,b)}
\epsilon_\gamma. \]
For points $a$ and $b$ in ${\cal M}_C$, with $\mu(a)-\mu(b)=1$ mod
$l$,  define the boundary components by considering the one
dimensional component ${\cal M}^1(a,b)$ of the moduli space ${\cal
  M}(a,b)$, 
\[ <\partial a, b>=\epsilon(a,b), \]
where
\[ \epsilon(a,b)=\sum_{\gamma\in\hat{\cal M}^1(a,b)}
\epsilon_\gamma. \]
\label{bound1}
\end{defin}

The construction of the Seiberg--Witten Floer homology is the result
of the following theorem.

\begin{thm}
The boundary operators of definiton \ref{bound1}
\begin{equation}
\partial^H a=\sum_{b\mid \mu(a)-\mu(b)=1} \epsilon^H(a,b) b
\label{boundH}
\end{equation}
and
\begin{equation}
\partial a=\sum_{b\mid \mu(a)-\mu(b)=1~mod~l} \epsilon(a,b) b
\label{bound}
\end{equation}
satisfy $\partial^H\circ \partial^H=0$ and $\partial\circ
\partial=0$. Upon fixing the Morse index of a 
critical point, it is possible to construct chain complexes with 
\[ C^H_q =\{ b\in {\cal M}_C^H \mid \mu(b)=q \}, \]
\[ C_q =\{ b\in {\cal M}_C \mid \mu(b)=q~mod~l \}, \]
and with the boundary operators
described above. The $\ZZ$-graded Seiberg-Witten Floer homology is defined as
\[ SWH_q(Y):= H_q(C^H_*,\partial^H)=
\frac{Ker(\partial^H_{q-1})}{Im(\partial^H_q)}. \] 
The $\ZZ_l$-graded Seiberg-Witten Floer homology is the homology
\[ SWH^l_q(Y):= H_q(C_*,\partial). \]
\label{homology}
\end{thm}

The fact that the boundary square is zero can be shown following
the analogous argument for the Donaldson case given in \cite{Fl}.
In fact, since $<\partial^H a,b>$ is the algebraic sum over the zero
dimensional manifold $\hat{\cal M}^H(a,b)$, the matrix elements for
$\partial^H\circ\partial^H$ are given by 
\[ <\partial^H\partial^H a,c>=\sum_b <\partial^H a,b><\partial^H b,c>. \]
This is the algebraic sum of the points in the the zero dimensional manifold
\[ {\cal M}^2(a,c):=\cup_{b}\hat{\cal M}^H(a,b)\times\hat {\cal M}^H(b,c), \]
with the product orientation $\epsilon^2(A,B)=\epsilon(A)\epsilon(B)$.

Now the claim follows if we prove that
the manifold ${\cal M}^2(a,c)$ with the above orientation is the oriented 
boundary of the 1-dimensional manifold $\hat{\cal M}^H(a,c)$.
This depends on a gluing formula like the one proven in
\cite{Fl}. The gluing formula will be proven in the following
section. 

The same argument applies to prove that $\partial^2=0$, according to
the spectral formula \ref{spectralsum}.

The compactness of the moduli space of solutions of (\ref{extrem})
implies that there is just a finite number of critical points of the
functional $\tilde {\cal C}$, and therefore only finitely many of the Floer
homology groups $SWH^l_*(Y)$ are non-trivial and they are finitely
generated. In the case of the Floer groups $SWH_*(Y)$, if $H$ is
infinite, there are infinitely many non-trivial groups. Each one is
finitely generated, but, if the group $H$ is non-trival, the
$\ZZ$-graded complex has generators in 
infinitely many degrees. In fact we have $C_{k}\cong C_{k+l}$, so that the
same groups appear with an $l$-periodicity.  
Notice that if $H$ is trivial then ${\cal M}_C^H$ is just ${\cal
  M}_C$, the Floer homology $SWH^l_*(Y)$ is $\ZZ$-graded and finitely
generated and coincides with $SWH_*(Y)$.  
Consider the homomorphism $\xi: H \to \ZZ$ given by
$$ \xi(h)=\langle c_1(L)\cup h, [X] \rangle. $$
The map $\xi$ is injective since $H=H^1(X,\ZZ)/\tilde
H^1(X,\ZZ)$, where $\tilde H^1(X,\ZZ)=ker(\xi)$. In particular
this means that $H$ is either trivial or isomorphic to $\ZZ$ via
the map $\xi$.  

\subsection{Convergence and Gluing}

The analysis presented in this section is part
of a joint work with B.L. Wang and will appear in \cite{MaWa}.

We show that the moduli spaces $\hat{\cal M}(a,b)$ of
unparametrized flow lines have a compactification with boundary strata
consisting of trajectories that break through other critical orbits.  
We discuss convergence in the unparametrized moduli spaces $\hat{\cal
  M}(a,b)$ and a gluing formula that describes the boundary
strata. 

We give the following preliminary definition.
\begin{defin}
A smooth path $(A(t),\psi(t))$ in ${\cal A}$ is of finite energy if
the integral
\beq
\int_{-\infty}^{\infty} \| \nabla \tilde{\cal C}(A(t),\psi(t))
\|^2_{L^2} dt <\infty  
\label{finenergy}
\eeq
is finite.
\end{defin}

Notice that any solution of (\ref{3SW1P'}) and (\ref{3SW2P'}) with
asymptotic values $a$ and $b$ is of finite energy, in fact
in this case the total variation of the functional $\tilde{\cal C}$
along the path $(A(t),\psi(t))$ is finite and
(\ref{finenergy}) satisfies
\[ \int_{-\infty}^{\infty} \| \nabla \tilde{\cal C}
(A(t),\psi(t)) \|_{L^2}^2 dt\leq C (\tilde{\cal C}(a)-\tilde{\cal C}
(b)), \] 
because of the assumptions on the perturbation $q_{(\AAA,\Psi)}$.
Finite energy solutions of the flow equations have nice properties:
they necessarily decay to asymptotic values that are critical points
of $\tilde{\cal C}$ as we prove in Proposition \ref{finenergy2}. We
begin by introducing some analytic properties of the functional $\tilde{\cal
  C}$ (see also \cite{Fr}, \cite{MST}, \cite{Wa}). 

\begin{lem}
Let ${\cal M}_C$ be the moduli space of critical points of $\tilde{\cal
  C}$, with $\rho$ a sufficiently small perturbation.
For any $\epsilon >0$ there is a $\lambda >0$ such that if the
$L^2_1$-distance from a point $[A,\psi]$ of ${\cal B}$ to all the
points in ${\cal M}_C$ is at least $\epsilon$, then 
\[ \| \nabla\tilde{\cal C}(A,\psi) \|_{L^2} >\lambda. \]
\label{PalaisSmale}
\end{lem}

\noindent\underline{Proof:} For a sequence $[A_i,\psi_i]$ of
elements of ${\cal B}$ with a distance at least $\epsilon$ from all
the critical points, such that
\[ \| \nabla\tilde{\cal C}(A_i,\psi_i) \|_{L^2}\to 0, \]
as $i\to\infty$, we would have
\[ \| *F_{A_i}-\sigma(\psi_i,\psi_i)-i\rho \| + \|
\partial_{A_i}\psi_i \| \to 0. \]
Thus, there is a constant $C$ such that
\[ \int_Y |*F_{A_i}-\sigma(\psi_i,\psi_i)-i\rho |^2 +
|\partial_{A_i}\psi_i |^2 dv <C. \]
If the perturbation $\rho$ is sufficiently small, the Weizenb\"ock
implies that
\[ \int_Y |F_{A_i}|^2 + |\sigma(\psi_i,\psi_i)|^2
+\frac{\kappa}{2}|\psi_i|^2 + 2|\nabla_{A_i}\psi_i|^2 dv <C. \]
Thus we have a uniform bound on the norms $\| \psi_i \|_{L^4}$, $\|
F_{A_i} \|_{L^2}$, and $\|\nabla_{A_i}\psi_i\|_{L^2}$. An elliptic
estimate shows that there is a subsequnce that converges in the
$L^2_1$ norm at a solution of the critical point equations
(\ref{extremP}), and this contradicts the assumption.

\noindent QED

\begin{corol}
Let $(A(t),\psi(t))$ be a smooth finite energy solution of equations
(\ref{3SW1P'}) and (\ref{3SW2P'}) with a smooth perturbation $q$.
Then there exist critical 
points $a$ and $b$ of $\tilde{\cal C}$, such that the 
$\lim_{t\to\pm\infty}(A(t),\psi(t))$ are in the gauge classes of $a$
and $b$.
\label{finenergy2}
\end{corol}

\noindent\underline{Proof:}
The finite energy condition (\ref{finenergy}) implies that
\[ \| \nabla\tilde{\cal C}(A(t),\psi(t)) \|\to 0 \]
as $t\to\pm\infty$. 
The Palais-Smale condition of Lemma \ref{PalaisSmale} implies that
there exist $T$ large, such that
for $|t| > T$,  $(A(t), \psi(t))$ lies in a very small
$\epsilon$-neighbourhood of critical points of $\tilde{\cal C}$.

\noindent QED

Now we can analyse convergence in the moduli space $\hat{\cal
  M}(a,b)$, see \cite{CaWa}, \cite{GWa}. 

\begin{thm}
The space $\hat{\cal M}(a,b)$ is precompact. Namely, any
sequence $x_i$ of elements in $\hat{\cal M}(a,b)$ has a
subsequence of smooth representatives that converges with all
derivatives to a solution $x$ of the flow equations which
lies in some $\hat{\cal M}(c,d)$ with
$\mu(a)>\mu(c)>\mu(d)>\mu(b)$. 
\label{seqcompact}
\end{thm}

\noindent\underline{Proof:}  
We choose a lift of the elements of $\hat{\cal M}(a, b)$ to
${\cal M}(a, b)$ such that the
gradient flow $[A(t), \psi(t)] $ has equal energy on
$(-\infty, 0]$ and on $[0, \infty)$:
\begin{equation}
\int_{-\infty}^0 \|\nabla \tilde{\cal C} ([A(t), \psi(t)] )\|_{L^2(Y)} dt
=\int_0^\infty\|\nabla \tilde{\cal C} ([A(t), \psi(t)] )\|_{L^2(Y)} dt.
\label{equal:energy}
\end{equation}

This lift is unique,
since in the family of gradient flows $\{ [\AAA, \Psi]^T, T\in\RR \}$
there is a unique element satisfying the equal energy condition
(\ref{equal:energy}).

Suppose $x_i=(\AAA_i, \Psi_i) \in {\cal A}_{k, \delta}(a, b)$
($i=1, 2, \cdots, \infty$) is 
a sequence of solutions to the equations (\ref{4SW1P'}) and
(\ref{4SW2P'}) which are represented in a temporal gauge by
the ``equal energy'' lifts $(A_i(t), \psi_i(t))$ of the sequence $x_i$ in
$\hat{\cal M}(a, b)$. 
This implies that the $(\AAA_i, \Psi_i)$ have a uniformly finite
energy $E$. 
By the Palais-Smale condition (Lemma \ref{PalaisSmale}),
we can find $T >> 1$ (choose $T > E/\lambda$ where
$\lambda$ is the constant appearing in Lemma \ref{PalaisSmale}), such that
for $|t| > T$, the $[A_i(t), \psi_i(t)]$ lie in a very small
$\epsilon$-neighbourhood of $a$ or $b$.

Therefore, the $[A_i(t), \psi_i(t)]$ have a uniform exponential decay 
over $(-\infty, -T]$ and $[T, \infty]$.
On $Y\times [-T-1, T+ 1]$ we proceed with an
argument that is analogous to the usual proof of the compactnes
for the Seiberg-Witten moduli space on a compact 4-manifold \cite{KM},
\cite{Mo}. Namely, we have the following.

A uniform $L^2$ bound on the spinors $\Psi_i$ follows from the
Weizenb\"ock formula 
\[ D_{\AAA}D_{\AAA}^*\Psi=\nabla_{\AAA}\nabla_{\AAA}\Psi +\frac{1}{2}
F_{\AAA}^+\Psi +\frac{\kappa}{4}\Psi, \]
where $\kappa$ is the scalar curvature. We have
\[ \| \Psi_i \|\leq \max_{Y\times [-T-1, T+ 1]} (-\kappa, 0) +2\|
P_{(\AAA_i,\Psi_i)}\|. \]
The perturbations $P_{(\AAA_i,\Psi_i)}$ are bounded uniformly
with respect to $(\AAA_i,\Psi_i)$ by assumption, \ref{calP}.
This also give a uniform bound on $\| \nabla_A \Psi_i \|_{L^2}$, as in
\cite{Mo} Lemma 5.1.7(?). The presence of the perturbation does not
affect the estimates, because of the boundedness of $P_{(\AAA_i,\Psi_i)}$.

We have the following gauge fixing condition: up to gauge
transformations $\lambda_i$ in the 
group ${\cal G}_{k+1}(Y\times [-T-1, T+ 1])$, it is possible to make
$\AAA_i-\AAA_0$ into a sequence of co-closed 1-forms, with the property
that 
\[ \| \AAA_i-\AAA_0\|^2_{L^2_k}\leq C\| F_{\AAA_i}^+ \|_{L^2_{k-1}}+K. \]
This is proven in Lemma 5.3.1 of \cite{Mo}.

{}From the curvature equation we have 
\[ \| F_{\AAA_i}^+\|\leq \| \Psi_i\cdot\bar\Psi_i \| + \| i\mu
+P_{(\AAA_i,\Psi_i)} \|, \]
which gives a uniform $L^2$-bound on $\| F_{\AAA_i}^+\|$. By the gauge
fixing condition this provides an $L^2_2$-bound on the $\AAA_i-\AAA_0$.

The elliptic estimates
\[ \|\AAA_i-\AAA_0\|_{L^2_k} \leq C\left(\|(d^*
+d^+)(\AAA_i-\AAA_0)\|_{L^2_{k-1}}+  
\|\AAA_i-\AAA_0\|_{L^2_{k-1}}\right) \]
and 
\[ \| \Psi_i \|_{L^2_k}\leq C \left( \| \nabla_A \Psi_i \|_{L^2_{k-1}}
  +\| \Psi_i \|_{L^2_{k-1}} \right) \]
provide a uniform bound on the $L^2_k$-norms of the $(\AAA_i,\Psi_i)$
on $Y\times [-T-1, T+ 1]$.  
By the Sobolev embeddings, this implies that on $Y\times
[-T-1, T+ 1]$ there is a subsequence $(\AAA_{i'},\Psi_{i'})$ that converges
uniformly with all derivatives.

Thus we can show that $(\AAA_i, \Psi_i)$ has
subsequence that converges to  a solution of (\ref{4SW1P'}) and
(\ref{4SW2P'}) on $Y\times R$. In fact,
by the uniform exponential decay, there
is  a subsequence converging
strongly on $Y \times ((-\infty, -T] \cup [T , \infty))$.
On  $Y\times [-T-1, T+ 1]$, as we have seen, after passing to a
further subsequence, 
there exist gauge transformations $u_i \in {\cal G}^0_{k+1}(Y\times
[-T-1, T+ 1])$ 
such that the transformed solutions converge strongly on
$Y\times [-T-1, T+ 1]$. We need to
merge $\{ u_i \}$ on
the overlap $K = Y \times ([-T-1, -T] \cup [T, T+1])$. This can be done
by choosing a cut-off function $c$ equal to 1 on $Y\times [-T, T]$ and
to 0 on $Y\times ((-\infty, -T-1] \cup [T+1 , \infty))$.
Over  $K$, there exists a subsequence of $\{ u_i \}$ converging strongly
to a gauge transformation $u$. For a sufficiently large $N$ and for all $i>N$,
 we have the $C^0$-bound $|u_i - u_N| < 1/2$. Then for all $i>N$,
$u_i = u_N exp (2\pi i \theta)$ for a unique fucntion $\theta_i$ on $K$
satisfying $|\theta_i| < 1/2$. Now define
gauge transformations $\{ v_i\}$ on $Y\times R$ by
\[
v_i = \left\{
\begin{array}{lll}
u_i &\qquad &\hbox{on $Y \times [-T, T]$,}\\[2mm]
 u_N exp (2\pi i c \theta_i)
     &\qquad& \hbox{on $Y \times ([-T-1, -T] \cup [T , T+1])$,}
\\[2mm]
u_N &\qquad& \hbox{on $Y \times ((-\infty, -T-1] \cup [T+1 , \infty))$.}
\end{array}\right.
\]
Then the sequence $v_i (\AAA_i, \Psi_i)$ converges strongly on $Y\times \RR$.
We denote with $(\AAA_\infty, \Psi_\infty)$ the limit. This satisfies
equations (\ref{4SW1P'}) and (\ref{4SW2P'}) on $Y\times \RR$ and is of
finite energy. Thus by Corollary \ref{finenergy2} it determines an
element of $\hat{\cal M}(c,d)$. For dimensional reasons we must
have $\mu(a)>\mu(c)>\mu(d)>\mu(b)$.

\noindent QED

The same argument applies to the moduli spaces $\hat{\cal M}^H(a,b)$,
see \cite{GWa}.

Notice that there can be at most finitely many distinct orbits $c$,
$d$ in ${\cal B}^H$ that appear as endpoints of limits of
sequences of solutions in $\hat{\cal M}^H(a,b)$. This happens because we
impose the condition $\int_Y c_1(L)\wedge [\lambda]=0$
for the gauge tranformations on $Y$. In fact, with this condition 
no two gauge equivalent points $b$ and $\lambda\cdot b$ can satisfy
$\mu(a)-\mu(b)=\mu(a)-\mu(\lambda\cdot b)$. 
Thus there can only be finitely many possible $\hat{\cal M}^H(c,d)$
in the boundary of $\hat{\cal M}^H(a,b)$. 

If one considers the identity
component in the group of gauge transformations on $Y$ and the
critical orbits in the space ${\cal B}^0$ of connections and sections
modulo the action of the connected component of the based gauge group, then 
there may be infinitely many orbits $\lambda\cdot b$ that have the
same index. In fact, suppose the group 
\[ H=\{ h\in H^1(Y,\ZZ) | \langle c_1(L)\cup h, [Y]\rangle =0 \]
is infinite. (This is obviously the case for instance when $c_1(L)=0$
rationally.) Then there are infinitely many distinct critical orbits
$\lambda\cdot b$ in ${\cal B}^0$ with
$\mu(a)-\mu(b)=\mu(a)-\mu(\lambda\cdot b)$. Thus, all the
moduli spaces $\hat{\cal M}(a, \lambda\cdot b)$ can appear in the
boundary of $\hat{\cal M}(a,b)$. In this case $\hat{\cal M}(a,b)$
does not have a nice compactification.

Theorem \ref{seqcompact} proves that lower dimensional moduli spaces
appear naturally in 
the compactification of the spaces $\hat{\cal M}^H(a,b)$. In the rest
of this section we describe a gluing formula, thus proving that the
boundary strata consist precisely of broken trajectories that live in
lower dimensional moduli spaces.

\begin{thm}
Suppose given $a$, $b$ and $c$ in ${\cal M}^H_C$  with
$\mu(a)>\mu(b)>\mu(c)$.  
Assume that $b$ is irreducible. Then, given a compact set 
\[ K\subset {\cal M}^H(a,b) \times {\cal M}^H(b,c), \]
there are a lower bound $T(K)>0$ and a smooth map 
\[ \#: K\times [T(K),\infty)\to {\cal M}^H(a,c) \]
\[ ((\AAA_1,\Psi_1),(\AAA_2,\Psi_2),T)\mapsto
(\AAA_1\#_T\AAA_2,\Psi_1\#_T\Psi_2), \] 
such that $\#_T$ is an embedding for all $T>T(K)$.
The gluing map $\#$ induces a smooth embedding
\[
\hat\#:  \hat K\times [T(K),\infty) \to \hat{\cal M}^H(a,c),
\]
where
\[  \hat K\subset \hat{\cal M}^H(a,b) \times \hat{\cal
  M}^H(b,c). \]
\label{equivgluing}
\end{thm}

This gives the compactification of the space
$\hat{\cal M}^H(a,b)$. The same result applies to the moduli space
$\hat{\cal M}(a,b)$ and gives an analogous compactification.

\begin{corol}
For a generic choice of the metric and of the perturbation, $\hat{\cal
  M}^H(a,b)$ has a compactification obtained by adding boundary
strata of the form
\[
\bigcup_{c_1, \cdots c_k} \hat{\cal M}^H(a, c_1) \times 
\hat{\cal M}^H({c_1},{c_2})\times \cdots \times
\hat{\cal M}^H({c_k},b).
\]
Here the union is over all possible sequences of the critical points
$c_1, \cdots c_k$ with decreasing indices.  
\label{strata}
\end{corol}

\noindent\underline{Proof of Theorem \ref{equivgluing}:}
The proof consists of several steps. We first define a pre-gluing map
$\#_T^0$ which provides an approximate solution, and then we prove
that this can be perturbed to an actual solution. We follow the
similar argument presented in \cite{Schwartz}.

The pre-gluing map is defined via the following construction. 
Let $x_1(t)=(A_1(t),\psi_1(t))$ and $x_2(t)=(A_2(t),\psi_2(t))$ be
elements in the moduli spaces 
${\cal M}^H(a,b)$ and ${\cal M}^H(b,c)$ respectively. 

We can write $x_1(t)=b+(\alpha_1(t),\phi_1(t))$ and
$x_2(-t)=b+(\alpha_2(t),\phi_2(t))$, where the elements
$(\alpha_i,\phi_i)$ have exponentially decaying behaviour, as in
Proposition \ref{decay}.
We construct an approximate solution $x(t)=(A_1\#_T^0 A_2 (t),
\psi_1\#_T^0 \psi_2 (t))$ as in \cite{Ma}, \cite{Wa},
of the form
\[ x(t)=\left\{ \begin{array}{lr} (A_1(t+2T),\psi_1(t+2T))&t\leq -1\\
b+\rho^-(t)(\alpha_1(t+2T),\phi_1(t+2T))+& \\
\rho^+(t)(\alpha_2(t-2T),\phi_2(t-2T))& -1\leq t\leq 1 \\ 
(A_2(t-2T),\psi_2(t-2T))&t\geq 1.
\end{array}\right. \]
Here $\rho^\pm(t)$ are smooth cutoff functions with bounded
derivative, such that $\rho^-(t)$ is equal
to one for $t\leq -1$ and to zero for $t\leq 0$ and $\rho^+(t)$ is
equal to zero for $t\leq 0$ and to one for $t\geq 1$. 

Consider the Hilbert bundles ${\cal T}_1$ and ${\cal T}_0$ that are
defined respectively as pullbacks via the map $\#_T^0$ of the
$L^2_{1,\delta}$ and of the $L^2_{0,\delta}$ tangent bundles of ${\cal
  B}^H(a,c)$, on the base space $K\times [T_0,\infty)$.

The flow 
\[
\frac{d}{dt}\psi_1 \#^0_T \psi_2=-\partial_{A_1 \#^0_T A_2}\psi_1\#^0_T\psi_2,
\]
\[ \frac{d}{dt}A_1\#^0_T A_2
=\sigma(\psi_1\#^0_T\psi_2,\psi_1\#^0_T\psi_2) -*F_{A_1\#^0_T A_2}
+2i\rho +2q_{(\AAA_1\#^0_T\AAA_2,\Psi_1\#^0_T\Psi_2)} \]
defines the fibre restriction of a bundle map from ${\cal T}_1$ to
${\cal T}_2$ defined on a neighbourhood of the zero section in ${\cal
  T}_1$. The linearization ${\cal L}_x$ at the approximate solution $x$ is the
fibre derivative of the above bundle map.

Since the linearizations ${\cal L}_{\AAA_1,\Psi_1}$ and ${\cal
    L}_{\AAA_2,\Psi_2}$ are surjective, then 
\[ {\cal K}=\bigcup_{K\times [T_0,\infty)} Ker({\cal
  L}_{\AAA_1,\Psi_1})\times Ker({\cal L}_{\AAA_2,\Psi_2}) \]
is a subbundle of ${\cal T}_1$. Thus we can define a space
${\cal T}^\perp_\chi$ for $\chi\in K\times [T_0,\infty)$ given by all
elements of ${\cal T}_1$ that are orthogonal to the image of
$Ker({\cal L}_{\AAA_1,\Psi_1})\times Ker({\cal L}_{\AAA_2,\Psi_2})$
    under the linearization $L_{\#}$ of the pre-gluing $\#_T^0$.

\begin{prop}
There exist a bound $T(K)\geq T_0$ such that, for all $T\geq T(K)$ and
for all broken trajectories 
\[ \left((A_1(t),\psi_1(t)),(A_2(t),\psi_2(t))\right)\in K, \] 
the Fredholm operator ${\cal L}_x$ 
\[ {\cal L}_x: {{\cal T}_1}_x\to {{\cal T}_0}_x \]
is surjective, where $x(t)$ is the approximate solution.
Moreover, composition of the pre-gluing map $\#^0_T$ with the
orthogonal projection on $Ker {\cal L}_x$ gives an isomorphism
\[ Ker({\cal L}_{\AAA_1,\Psi_1})\times Ker({\cal L}_{\AAA_2,\Psi_2})
    \stackrel{\cong}{\to} Ker{\cal L}_x. \]
\label{approxsurjective}
\end{prop}

\noindent\underline{Proof of Proposition \ref{approxsurjective}:} We
know that ${\cal L}_x$ is Fredholm of index $\mu(a)-\mu(c)$. We
also know that $\dim Ker({\cal L}_{x_1})=\mu(a)-\mu(b)$
  and $\dim Ker({\cal L}_{x_2})=\mu(b)-\mu(c)$.

We need to show that for any pair $((\AAA_1,\Psi_1),(\AAA_2,\Psi_2))$
there is a bound $T_0=T(x_1,x_2)$ such that ${\cal
  L}_x$ is surjective for $T\geq T_0$. The compactness of $K$ will
ensure that there is a uniform such bound $T(K)$.

It is therefore enough to prove the following crucial step.

\begin{lem}
There exist $T_0$ and a constant $C>0$ such that 
\[ \| {\cal L}_x \xi \|_{L^2_{0,\delta}}\geq C \|\xi
\|_{L^2_{1,\delta}}. \] 
for all $T\geq T_0$ and $\xi$ in ${\cal T}^\perp_\chi$, where in our
notation $x=\#_T^0(\chi)$.  
\label{surj2}
\end{lem}

\noindent\underline{Proof of Lemma \ref{surj2}:}
Suppose there are sequences $T_k\to\infty$ and $\xi_k\in{\cal
  T}^\perp_\chi$ such that $\| \xi_k \|=1$ and $\| {\cal L}_x \xi_k \|\to
0$. 

We first show that the supports of the $\xi_k$ become more and more
concentrated at the asymptotic ends as $k\to\infty$.
We consider the operator ${\cal L}_b=\frac{\partial}{\partial
  t}+L_b$. 
If $\zeta:\RR\to[0,1]$ is a smooth function which is equal to 1 on
$[-1/2,1/2]$ and equal to zero outside $(-1,1)$, let
$\zeta_k(t)=\zeta(\frac{t}{2T_k})$. Then we have
\[ \| {\cal L}_b ~\zeta_k\xi_k \|\leq \| \zeta_k^\prime \xi_k \| + \| \zeta_k
{\cal L}_b \xi_k \| \leq \]
\[ \leq \frac{1}{2T_k} \max |\zeta^\prime | +\|({\cal L}_x-{\cal L}_b)
\xi_k \| +\| {\cal L}_x \xi_k \|\leq  \]
\[ \leq \frac{1}{T_k} \max |\zeta^\prime |+ \sup_{t\in[-T_k,T-k]}\|
L_{x(t)}-L_b\|~ \|\xi_k\| +\| {\cal L}_x \xi_k \|. \]
Thus, $\| {\cal L}_b ~\zeta_k\xi_k \|\to 0$ as $k\to\infty$. In fact,
the term $\sup_{t\in[-2T_k,2T-k]}\| L_{x(t)}-L_b\|$ is bounded by
\[ \sup_{t\in[-1,1]}\| L_{x(t)}-L_b\| + \sup_{t\in[-2T_k,-1]}\|
L_{x_1(t+2T_k)}-L_b\| +\sup_{t\in[1,2T_k]}\| L_{x_2(t-2T_k)}-L_b\|. \]
All these terms tend to zero because of the exponential decay to the
critical point $b$ of the trajectories $(A_1(t),\psi_1(t))$ and
$(A_2(t),\psi_2(t))$. The operator ${\cal L}_b$ is an isomorphism 
between the spaces $L^2_{1,\delta}(\RR,T({\cal B}^H))$ and
$L^2_{0,\delta}(\RR,T({\cal B}^H))$, hence we have $\xi_k\to 0$ in the
$L^2_{1}$ norm over $Y\times [-2T_k,2T_k]$. 

This result allows us to rephrase the convergence condition $\| {\cal
  L}_x \xi_k \|\to 0$ in terms of the Fredholm operators
${\cal L}_{x_1}$ and ${\cal L}_{x_2}$: 
\[ \| {\cal L}_{x_1} (\rho^-_{1-T_k}\xi_k^{-T_k}) \|\leq \|{\rho^\prime}^-
\xi_k\| + \| \rho^- {\cal L}_x \xi_k \| \leq \]
\[ \leq  C \|\xi_k\|_{Y\times [-1,1]} + \| {\cal L}_x \xi_k \| \to
0, \]
where $\rho^-_{1-T_k}(t)=\rho(t+1-T_k)$ and $\xi_k^{-T_k}(t)=\xi_k(t-T_k)$.
This implies $\rho^-_{1-T_k}\xi_k^{-T_k}\to v$ where $v\in Ker({\cal
  L}_{x_1})$, since ${\cal L}_{x_1}$ is a Fredholm operator. Thus, $\|
\rho^-_1 \xi_k- v^{T_k}\|\to 0$. Similary we obtain an element $u$ in
$Ker({\cal L}_{x_2})$ such that $\| \rho^+_{-1}\xi_k -u^{-T_k}\|\to 0$.

We now use these estimates to derive a contradiction with the
assumption that $\| \xi_k \|=1$ and $\xi_k\in{\cal T}^\perp_\chi$.
We have
\[ 1=\lim_k \|\xi_k\|=\lim_k \langle \rho^-_1 \xi_k,\xi_k\rangle + \langle
\rho^+_{-1}\xi_k,\xi_k\rangle, \]
since the remaining term satisfies 
\[ \langle (1-\rho^-_1-\rho^+_{-1})\xi_k,\xi_k\rangle=0 \]
for large $k$ because $(1-\rho^-_1-\rho^+_{-1})$ is supported in
$[-2,2]$. 
Thus the equality can be rewritten as
\[ 1= \lim_k \langle \rho^- v,\xi_k \rangle +\lim_k \langle \rho^+
u,\xi_k \rangle =
\lim_k \langle L_{\#}(u,v),\xi_k \rangle=0. \]
The last equality holds since, by construction, $\xi_k\in {\cal
  T}^\perp$ is orthogonal to the image $L_{\#}$ 
of the linearization of the pre-gluing map. 

This completes the proof of Lemma \ref{surj2}. Now the rest of
Proposition \ref{approxsurjective} follows, since we obtain
\[ \dim Ker({\cal L}_x )= \dim Ker ({\cal L}_{\AAA_1,\Psi_1}) + \dim
Ker({\cal L}_{\AAA_2,\Psi_2}). \]

\noindent QED

This means, as we are going to discuss, that ${\cal T}^\perp$ is the
normal bundle for the gluing construction. 

Now we want to define the actual gluing map $\#$ that provides a
solution of the flow equations in ${\cal M}^H(a,c)$.
This means that we want to obtain a section $\sigma$ of ${\cal T}_1$ such
that the image under the bundle homomorphism given by the flow
equation is zero in ${\cal T}_0$.
Moreover, we want this element $\sigma(\AAA_1,\Psi_1,\AAA_2,\Psi_2,T)$ to
converge to zero sufficiently rapidly as $T\to\infty$, so that the
glued solution will converge to the broken trajectory in the limit
$T\to\infty$. 

The result is obtained as a fixed point theorem in Banach spaces. 
Consider the right inverse map of ${\cal L}$ restricted to ${\cal T}^\perp$, 
\[ G: {\cal T}_0\to {\cal T}^\perp. \]
There is a $T(K)$ and a constant $C>0$ such that
\[ \|G_\chi\xi\|\leq C \|\xi \| \]
for $\chi\in K\times [T(K),\infty)$. Some care is needed in
obtaining the uniformity of the constant $C$ with respect to
$\chi$. We refer to \cite{Schwartz} for further details.

We aim at using the contraction principle. Namely, suppose given a
smooth map $f:E\to F$ between Banach spaces of the form
\[ f(x)=f(0)+Df(0)x+N(x), \]
with $Ker(Df(0))$ finite dimensional, with a right inverse
$Df(0)\circ G=Id_F$, and with the nonlinear part $N(x)$ satisfying the
estimate 
\beq
\label{nonlinestimate}
\| GN(x)-GN(y) \|\leq C(\|x\|+\|y\| )\|x-y\| 
\eeq
for some constant $C>0$ and $x$ and $y$ in a small neighbourhood
$B_{\epsilon(C)}(0)$. Then, with the initial condition $\| G(f(0))
\|\leq \epsilon/2$, there is 
a unique zero $x_0$ of the map $f$ in  $B_\epsilon(0)\cap G(F)$. This
satisfies $\| x_0 \|\leq \epsilon$.

The map $f$ is given in our case by the flow equation, viewed as a
bundle homomorphism ${\cal T}_1\mapsto {\cal T}_0$. We write $f$ as a
sum of a linear and a non-linear term, where the linear term is ${\cal
  L}$ and the nonlinear term is
\[ N_{(A(t),\psi(t))}(\alpha,\phi)=(\sigma(\phi,\phi)+{\cal
  N}q_{(A(t),\psi(t))}(\alpha,\phi),\alpha\cdot\phi). \]
Here we write the perturbation $q$
of equation (\ref{3SW2P'}) as sum of a linear and a non-linear term, 
$2q={\cal D}q+{\cal N}q$. 

It is clear that the quadratic form
$(\sigma(\phi,\phi),\alpha\cdot\phi)$ satisfies an
estimate of the form (\ref{nonlinestimate}). The perturbation term
also satisfies a similar estimate for 
large enough $T$ because of the assumptions on the perturbation space
${\cal P}$. This implies the estimate (\ref{nonlinestimate}) for $G\circ N$.

We have to verify the initial condition. This is provided by the
exponential decay. In fact, we have
\[ \| f(A_1\#_T^0 A_2,\psi_1\#_T^0\psi_2)\|\leq
C(\|(\alpha_1,\phi_1)^{2T}\|_{Y\times [-1,0]}
+\|(\alpha_2,\phi_2)^{-2T}\|_{Y\times [0,1]}). \]
The exponential decay of $(A_1(t),\psi_1(t))$ and $(A_2(t),\psi_2(t))$
to the endpoints implies a decay
\beq
\| f(A_1\#_T^0 A_2,\psi_1\#_T^0\psi_2)\|\leq C
e^{-\delta T} 
\label{gluedecay}
\eeq
For all $T\geq T_0$. The constant $C$ and the lower bound $T_0$ can be
chosen uniformly due to the compactness of $K$.

This provides the existence of a unique correction term 
\[ \sigma(\AAA_1,\Psi_1,\AAA_2,\Psi_2,T)\in B_\epsilon(0)\cap {\cal
  T}^\perp \] 
satisfying $f(\sigma)=0$. The implicit function theorem ensures that
$\sigma$ is smooth. The exponential decay (\ref{gluedecay})
ensures an analogous decay for $\sigma$, hence the glued trajectory
approaches the broken trajectory when $T$ is very large. The gluing
map is given by
\[ (A_1\#_T A_2,\psi_1\#_T\psi_2)=(A_1\#_T^0 A_2,\psi_1\#_T^0\psi_2) +
\sigma(\AAA_1,\Psi_1,\AAA_2,\Psi_2,T). \]

\noindent QED

\section{A Casson-like Invariant}

One can define an invariant by taking the Euler characteristic of the 
$SWH_*$ groups.  I have learned from B.L. Wang that the
construction of this invariant
have been worked out independently in \cite{ZW2}, which is a more
recent version of \cite{ZW} where the invariant was
introduced  in terms of the partition function of a Topological
Quantum Field Theory. 
 
In the case of Donaldson-Floer theory, there is a nice identification
of flat connections (critical points of the Chern-Simons functional)
with representations of the fundamental group in $SU(2)$. Therefore it
is possible to describe the invariant in purely geometric terms, and
it turns out to be the Casson invariant of homology 3-spheres, as
shown in \cite{T1}.

In our case there is no immediate geometric interpretation of the
solutions of (\ref{extrem}) in non-gauge theoretic terms.

However some of the results of \cite{T1} can be carried over to the
present case, with some minor modifications in the arguments. We shall
present these results in the next section.

The following result makes it, in principle, easier to compute this invariant
in the case when ${\cal M}_C$, the moduli space of critical points, is zero
dimensional.

\begin{thm}
The invariant $\chi(SWH_*(Y))$
is just the sum over points in ${\cal M}_C$ of the signs given
by the orientations defined in section 2.4,
\[ \chi(SWH_*(Y))= \sum_{a\in{\cal M}_C} \epsilon_a. \]
\label{inv0dim}
\end{thm}

\noindent\underline{Proof:}
Since the dimension of the set of critical points is zero, the
operator $L$ has trivial kernel. This implies $Ker(T)=0$ at the
critical points. Given a path $\gamma$ between two critical points,
consider the set $\Omega=\gamma\times \CC$.  Since $Ker(T_0)=Ker(T_1)=0$
(see \cite{APS}), the spectral flow of 
\[ T_t=T\mid_{(A(t),\psi(t))} \]
along $\gamma$ can be thought of as the
algebraic sum of the intersections in $\Omega$ of the set 
\[ {\cal S}=\cup_{t\in\gamma}Spec(T_t) \]
with the line $\{(t,0)\mid t\in\gamma \}$.
This counts the points where the discrete spectrum of the operator $T$
crosses zero, with the appropriate sign. Up to perturbations we can
make these crossings transverse.

We can express the same procedure in a different, yet equivalent, way.
Consider $T_t$ as a section over ${\cal B}$ of the bundle of index
zero Fredholm operators, $Fr_0$.
There is a first Stiefel-Whitney class $w_1$ in $H^1(Fr_0;\ZZ_2)$ that
classifies the determinant line bundle of $Fr_0$ (see \cite{T1}, \cite{Ko}).
 
The submanifold of codimension one in $Fr_0$ that represents
the class $w_1$ is given by Fredholm operators of nonempty kernel,
$Fr_0^1\subset Fr_0$. This submanifold can be thought of as the zero
set of a generic section of the determinant line bundle. 

Given a path $\gamma :I\to {\cal B}$ joining two points $a$ and $b$, 
the image of $\gamma$ composed
with a generic section $\sigma$ of $Fr_0$ will meet $Fr_0^1$ transversally.

We call this intersection number $\delta_\sigma (a,b)$. 
This counts the points in $Det(\sigma\circ\gamma)^{-1}(0)$ with the
orientation. If we take the
section $\sigma\circ\gamma$ to be $T_t$ what we get is exactly the
spectral flow between the critical points $a$ and $b$.
This number is defined modulo $l$, but the mod 2 spectral flow
does not depend on the choice of the path.
Now we want to show that this same intersection number measures the
change of orientation between the two points $a$ and $b$ of ${\cal
  M}_C$.

The orientation of the tangent space at a critical point is
given by a trivialization of $Det(L)$ or equivalently of $Det(T)$ at
that point.  

Up to a perturbation, we can assume that $Coker(\tilde T_t)$ is
trivial  along the path $\gamma$, see lemma \ref{dim1}. 
Therefore we can think of $Det(T_t)$
as $\Lambda^{top}Ker(T_t)$, which specifies the orientation of the
space $Ker(T_t)$. Thus we obtain that the change of orientation
between two critical points is measured by the mod 2 spectral flow.

This means that the sign attached to the point in the grading of the
Floer complex is the same as the sign that comes from the orientation
of ${\cal M}_C$.

\noindent QED

A definition of the Seiberg-Witten invariant of three-manifolds as the number
obtained by counting points in ${\cal M}_C$ with the orientation was given
in \cite{Au2}. 
   
The invariant $\chi(SWH_*(Y))$ vanishes for all $Y$ which admit a metric of
positive scalar curvature. This is a consequence of the Weitzenb\"ock
formula, \cite{W2}. However the invariant is non-trivial. In fact
consider the case of a 3-manifold $Y=\Sigma\times S^1$, with the
$Spin_c$ structure determined by the pullback of a line bundle $L$ of
degree $d$ on the surface $\Sigma$. As proved in \cite{MST} (see also
\cite{Do2} and \cite{Mu}), the Seiberg-Witten-Floer homology of $Y$ is
the ordinary homology of the symmetric product $s^r(\Sigma)$, with
$r=(2g-2)-d >0$. But the Euler characteristic of $s^r(\Sigma)$ is
\[ \chi(s^r(\Sigma))=\chi(SWH_*(\Sigma\times
S^1))=\left(\begin{array}{c} 2g-2 \\ r \end{array}\right). \]
Notice, however, that in this particular case the invariant is
computed with the unperturbed equations, which give rise to a positive
dimensional moduli space ${\cal M}_C$. This is a therefore an analogue
of a degenerate Morse theory. 

\subsection{Invariance}

We prove here that the Euler characteristic of the Floer homology 
is a topological invariant of the 3-manifold $Y$, if the manifold has
Betti number $b^1(Y)>1$. This means that the Euler characteristic of
the Floer homology is independent of
the choice of the metric and of the perturbation. We give an argument
which is similar to the one constructed for Donaldson theory in
\cite{DK}, pg. 140-146.

Suppose we are given two Riemannian metrics on $Y$, $g_0$, $g_1$. Consider a
path of metrics $g_t$, $t\in [0,1]$, connecting them.
Take the infinite dimensional manifold ${\cal B}\times [0,1]$.
We want to construct a ``universal moduli space'' for the equation
(\ref{extrem}). Consider the product bundle over  ${\cal B}\times
[0,1]$ with fibre $\Lambda^1(Y)\oplus \Gamma (S\otimes
L)$. We have a section given by
\begin{equation}
\label{univsec} 
s(A,\psi, t)=(\partial^t_A\psi, *_tF_A-2\rho_t -\sigma(\psi,\psi)). 
\end{equation}
Note that in  (\ref{univsec}) both the Dirac operator and
the Hodge $*$-operator depend on the metric, hence on the parameter $t$.
Here $\{ \rho_t \}$ is any family of 1-forms that are co-closed
with respect to $*_t$ and away from the wall, i.e. $*_t\rho_t\neq \pi c_1(L)$.
The universal moduli space is the zero set of the section $s$, ${\cal
  MU}=s^{-1}(0)$.

We need the following lemma. Recall that a map of Banach manifolds is
Fredholm if its linearization is a Fredholm operator.

\begin{lem}
The section $s$ is a Fredholm section. The operator that linearizes
(\ref{univsec}) is onto.
\label{fred}
\end{lem}

\noindent\underline{Proof:} The linearization of (\ref{univsec}) is given by 
\[ ds_{(A,\psi,\rho,t)} (\alpha,\phi,\eta,\epsilon)=\epsilon
\frac{\partial}{\partial t}s(A,\psi,\rho,t) +
\tilde T\mid_{(A,\psi,\rho,t)} (\alpha,\phi,\eta), \]
where $(\alpha,\phi)$ are coordinates in the tangent space of ${\cal
  B}$, $\eta$ is a 1-form, and $\epsilon\in \RR$.

The operator $\tilde T$ is the perturbed Hessian, 
\[ \tilde T\mid_{(A,\psi,\rho,t)} (\alpha,\phi,\eta)= -2\eta +
T\mid_{(A,\psi,\rho,t)}(\alpha,\phi). \]
We follow the same argument used in lemma \ref{dim1} to show that
$\tilde T$ is surjective for a generic choice of the perturbation.

\noindent QED

Now we can apply again the implicit function theorem on Banach manifolds,
and we get that $s$ is transverse to the zero section. Therefore 
the universal moduli space is a smooth manifold, provided
the perturbation $\rho_t$ is chosen in such a way that the moduli space
corresponding to each value of $t$ does not contain reducibles
($\psi\equiv 0$). It is sufficient to add the constraint that $*\rho_t\neq \pi
c_1(L)$. The dimension of the universal moduli space is 
$Ind(T)+1=1$. The proofs of the compactness and orientability of
the universal moduli space are analogous to the proofs given in section
2. 

The independence of the metric 
now follows from the fact that the moduli spaces corresponding to the
metrics $g_0$ and $g_1$ form the boundary of a compact
oriented 1-manifold, and the total oriented boundary of such a
manifold is zero.

\begin{rem}
\label{b1}
The condition that $*\rho_t\neq \pi c_1(L)$ is satisfied for a generic 
choice of the perturbation if the manifold $Y$ has $b^1(Y)>1$. If $Y$ has
$b^1(Y)=1$ then the condition is satisfied by sufficiently small 1-form
$\rho_t$, hence in this case the metric and the perturbation cannot be
both chosen arbitrarily. However, given $g_0$ and $g_1$ it is possible to
find sufficiently small perturbations such that the invariant does not
change along the path of metrics $g_t$. This has been pointed out
already in \cite{Au2}. 
\end{rem} 

\subsection{Heegaard Splittings}

In this section we extend to the case of Seiberg-Witten-Floer
homology a result of Taubes' \cite{T1} concerning the behaviour of
the Euler characteristic of Donaldson-Floer homology under a Heegaard
splitting of a homology 3-sphere.
This proves that the invariant $\chi(SWH_*(Y))$ behaves very much like
the Casson invariant. In fact it was conjectured by
Kronheimer \cite{KM2}, and recently proved by Chen \cite{Chen} and Lim
\cite{Lim} independently, that in the case of a homology sphere the Euler
characteristic of the Seiberg-Witten Floer
homology differs from the Casson invariant for a correction term that
depends on the index of the Dirac operator and on the signature of a
$Spin$-four manifold that bounds the homology sphere \cite{KM2}. 

We follow the same argument and the notation used in \cite{T1}. We
do not write out the analysis in full details in the proof
of the Heegaard splitting formula, since we prefer to concentrate on
the topological result and refer for the analysis to the case
of the instanton Floer homology and the Casson invariant as worked out
by Taubes \cite{T1}. 

Let us introduce some preliminary definitions.
Let ${\cal F}$ denote the 1-form on $\hat{\cal A}$ obtained as the 
differential of the functional $\tilde{\cal C}$ as in
(\ref{1forma}):
\[ {\cal F}\mid_{(A,\psi)}(\alpha,\phi)= -\int_Y \alpha\wedge (F_A
-2*\rho -*\sigma(\psi,\psi)) +\int_Y <\phi,\partial_A\psi>. \]

Consider the open sets $Y_1$, $Y_2$ and their intersection $Y_0$,
which is of the form $Y_0=\Sigma\times I$ for some interval $I$ and some
closed surface $\Sigma$ of genus $g$.

Consider the pullback of the spinor bundle under the inclusion maps
$i_l: Y_l\hookrightarrow Y$ and $j_l: Y_0\hookrightarrow Y_l$,
$l=1,2$. We shall indicate by ${\cal A}_l$, $l=0,1,2$, the space of
pullback connections and sections of the pullback bundles.
Consider the maps 
\[ {\cal A}\stackrel{i^*_1\times i^*_2}{\to}{\cal A}_1\times{\cal A}_2
\stackrel{j^*_1\times j^*_2}{\to} {\cal A}_0\times {\cal A}_0. \]
We shall denote by $\hat{\cal A}_l$ the set of irreducible pairs in
${\cal A}$ (which means $\psi$ not identically zero) that map to  irreducible
pairs in ${\cal A}_0$. Note that if a pair is a critical point of the
functional $\tilde{\cal C}$ on $Y$ then it maps to an irreducible
pair. In fact the section $\psi$ satisfies the equation
$\partial_A\psi=0$. Therefore if it vanishes on an open set it has
to be identically zero on all of $Y$.
Consider the induced maps on the quotient space $\hat{\cal B}_l$ of $\hat{\cal
  A}_l$ with respect to the action of the gauge group.
We define 1-forms ${\cal F}_l$ on the tangent spaces of
${\cal A}_l$ and induced forms on $\hat{\cal B}_l$ using the same
expression given in theorem \ref{hessian}. 
We shall use the notation ${\cal M}_l$ to denote the set ${\cal
  F}_l^{-1}(0)$ in $\hat{\cal B}_l$.

\begin{thm}
Consider a  Heegaard splitting of a closed oriented 3-manifold $Y$,
$Y=Y_1\cup_\Sigma Y_2$. The Euler characteristic 
\[ \chi=\chi(SWH_*(Y)) \]
of the Floer homology
is the intersection number of the manifolds $j^*_1{\cal M}_1$ and 
$j^*_2{\cal M}_2$ inside ${\cal M}_0$.
\label{heegsplit}
\end{thm}

\noindent\underline{Proof:}
To prove the theorem we need some preliminary steps.

The $L^2_k$-tangent spaces to the gauge orbits are
\[ {\cal T}_l =\{ (\alpha,\phi)\in \Lambda^1(Y_l)\oplus\Gamma(\tilde
S\otimes L\mid_{Y_l}), G^*(\alpha,\phi)=0, i^*_{\partial
  Y_l}(\alpha,\phi)=0 \}, \]
where we impose the vanishing condition on the boundary.
We shall also consider the Banach bundle ${\cal L}_l$ given by
\[ {\cal L}_l =\{ (\alpha,\phi)\in \Lambda^1(Y_l)\oplus\Gamma(\tilde
S\otimes L\mid_{Y_l}), G^*(\alpha,\phi)=0 \}, \]
with the $L^2_k$-norm and no boundary conditions.

\begin{lem}
In the sequence of maps
\[ {\cal A}\stackrel{I}{\to}{\cal A}_1\times{\cal A}_2
\stackrel{J}{\to} {\cal A}_0\times {\cal A}_0, \]
with $I=i^*_1\times i^*_2$ and $J=j^*_1\times j^*_2$, $I$ is an
embedding and $J$ is a submersion. Moreover, $Im(I)=J^{-1}(\Delta)$, i.e.
the image of $I$ is the inverse image of the diagonal under the map $J$.
Moreover
the same result holds for the induced sequence on the spaces $\hat{\cal
  B}_l$.
\label{l1}
\end{lem}

\noindent\underline{Proof:} We need to prove that
the induced maps have the required properties  on
the level of tangent spaces. This result can 
be obtained by the same proof as that given in \cite{T1}.
In fact the tangent maps induced by $I$ and $J$ can be written as
\[ I_*(\alpha,\phi)=(i_1^*(\alpha,\phi) +G(f_1),
i_2^*(\alpha,\phi) +G(f_2)) \]
and 
\[ J_*(\alpha_1,\phi_1,\alpha_2,\phi_2)=(j_1^*(\alpha_1,\phi_1)
+G(g_1), j_2^*(\alpha_2,\phi_2) +G(g_2)), \]
subject to the condition $G^*G(f_l)=G^*G(g_l)=0$ and the vanishing
condition on the boundary. The gauge maps are chosen so as to guarantee
that  $I_*$ and $J_*$ map tangent spaces to tangent spaces.

But the only solution to the equation $G^*Gf=0$ on $Y_l$ and 
$f=0$ on $\partial Y_l$ is the trivial one. Hence $I_*$ is
injective, $J_*$ is surjective and $Im (I_*)=Ker((j_1^*)_*-(j_2^*)_*)$.

\noindent QED

Define the forms $\nabla{\cal F}_l=\frac{d}{ds}{\cal
  F}_l((A,\psi)+s(\alpha,\phi))\mid_{s=0}$, as in (\ref{2forma}).
These quadratic forms define operators $T_l$ on the tangent space to
$\hat{\cal B}_l$ that coincide with the
Hessian of theorem \ref{hessian} when $(A,\psi)$ is a critical point
of $\tilde{\cal C}$, i.e. when ${\cal F}_l\mid_{(A,\psi)}\equiv 0$.

\begin{lem}
The operator $T_l$ is a Fredholm operator from 
\[ \Lambda^0(Y_l)\oplus\Lambda^1(Y_l)\oplus\Gamma(S \otimes
L\mid_{Y_l}), \]
completed in the norm $L^2_k$ with vanishing
conditions on the boundary, to the same space completed in the
$L^2_{k-1}$ norm. 
\end{lem}

\noindent\underline{Proof:} Up to a compact perturbation we have 
the operator $D$ of (\ref{D}) and the Dirac operator. The first is 
Fredholm from $\Lambda^0(Y_l)\oplus\Lambda^1(Y_l)$, completed in the
$L^2_k$ norm with vanishing conditions, to
$\Lambda^0(Y_l)\oplus\Lambda^1(Y_l)$, completed in the $L^2_{k-1}$ norm, and
the latter is Fredholm between the corresponding spaces of sections
$\Gamma(S \otimes L\mid_{Y_l})$.

\noindent QED

As a consequence of these lemmata we have the following result.

\begin{thm}
The sets ${\cal F}_l^{-1}(0)$ obtained from the
perturbed equations on $Y_l$ are embedded submanifolds ${\cal
  M}_{C,l}$ of $\hat{\cal B}_l$. 
Moreover the intersection of ${\cal M}_{C,1}$ and ${\cal
  M}_{C,2}$ in ${\cal M}_{C,0}$ is transverse.  
\end{thm}

\noindent\underline{Proof:} The first assertion of the theorem is an
application of the implicit 
function theorem for Banach manifolds and Fredholm operators. It is
sufficient to prove that $\dim Ker(T_l)=Ind(T_l)$.
Since we are assuming vanishing conditions on the boundary, surjectivity
of $T_l$ for a suitable choice of the perturbation follows
from the argument used in lemma \ref{dim1}.

For the trasversality property, we can prove the following claim.
The intersection is
transverse at a point $(A,\psi)\in {\cal M}_C$ if and only if the
kernel of $T$ at that point is trivial. 
But at a critical point this follows if we have that $Ker(L)$ is 
trivial. We proved in lemma \ref{dim1} that this is the case for 
a generic perturbation.

In order to prove the claim we follow \cite{T1} and introduce the operator  
\[ H\mid_{(A,\psi)}: {\cal T}_1\oplus{\cal T}_2\oplus{\cal L}_0\to
{\cal L}_1\oplus{\cal L}_2\oplus{\cal T}_0, \]
\[ H\mid_{(A,\psi)}(\alpha_l,\phi_l)= (T_1 (\alpha_1,\phi_1),T_2 
(\alpha_2,\phi_2), \]
\[ (j^*_1)_*(\alpha_1,\phi_1)-(j^*_2)_*(\alpha_2,\phi_2)) -T_0^*
(\alpha_0,\phi_0) ). \] 

The kernel of the operator $H$ measures the lack of transversality of
the intersection. In fact $Ker( H\mid_{(A,\psi)})$ is given by the pairs
$(\alpha_1,\phi_1)\in Ker(T_1)$ and
$(\alpha_2,\phi_2)\in Ker(T_2)$ such that the vector
\[ (j^*_1)_*(\alpha_1,\phi_1)-(j^*_2)_*(\alpha_2,\phi_2) \]
is orthogonal to
$Ker(T_0)$. This means exactly that the tangent spaces of
${\cal M}_{C,1}$ and ${\cal M}_{C,2}$ meet in a non
trivial subspace in order to sum up to the dimension of the
tangent space of ${\cal M}_{C,0}$.

Now suppose $(A,\psi)$ is a point in ${\cal M}_C$ and
$(\alpha_l,\phi_l)$ is in the kernel of $T_l$ for $l=1,2$.
It is not hard to check that 
\[ (j^*_1)_*(\alpha_1,\phi_1)-(j^*_2)_*(\alpha_2,\phi_2) \]
is in $Ker(T_0)$. 
So elements in the kernel of  $H\mid_{(A,\psi)}$ have 
\[ (j^*_1)_*(\alpha_1,\phi_1)-(j^*_2)_*(\alpha_2,\phi_2)=0. \]
This determines an element
$(\alpha,\phi)$ in $Ker(T)$. Thus at a critical point the
intersection is transverse if and only if $Ker(T)$ is
trivial.

\noindent QED

We can see with an index computation that the dimensions of the
intersecting manifolds match properly.

\begin{lem}
The sum of the indices of the operators $T_l$ over the
manifolds $Y_l$ with $l=1,2$ gives exactly the index of $T_0$ on $Y_0$.
\end{lem}

\noindent\underline{Proof:} The metric on $Y$ is such that it is a
cylinder on $Y_0$ and in a neighbourhood of the boundary on $Y_l$,
$l=1,2$. Thus on $\Sigma \times I$  we have
$Ind(D_0)=-\chi(\Sigma)=2g-2$ where $g$ is the genus of $\Sigma$.  
On $Y_l$ with $l=1,2$ we have $H^1(Y_l;\ZZ)=\ZZ^g$ and
$H^1(Y_l, \partial Y_l;\ZZ)=0$ since $Y_l$ is a handlebody.
Therefore $H^2(Y_l,\ZZ)=0$ and $g-1=-\chi(Y_l)=Ind(D_l)$.
The index of the Dirac operator is just the Atiyah-Patodi-Singer boundary term
on $\partial Y_l$, $l=1,2$. These sum up to give exactly the boundary
term on $\partial Y_0$. 

\noindent QED

From the previous results we have the following.

\begin{lem}
Up to a choice of the orientation, we can write the set of
critical points as
\[ {\cal F}^{-1}(0)= I^{-1}(({\cal F}_1^{-1}(0)\times{\cal
  F}_2^{-1}(0))\cap J^{-1}(\Delta)). \]
\label{points}
\end{lem}

The last step of the proof is to show that the sign difference of two
points in the oriented intersection is given exactly by the spectral
flow that defines the relative Morse index in the grading of
the Floer groups, up to an overall sign that comes from fixing the
Morse index of one particular solution.

The argument is analogous to the proof given in \cite{T1} adapted to the
present case along the line of the proof of theorem \ref{inv0dim}.

\noindent QED

\vspace{1cm}

\noindent Matilde Marcolli

\noindent Department of Mathematics

\noindent The University of Chicago

\noindent Chicago IL 60637 USA

\noindent matilde\@@math.uchicago.edu

\end{document}